

\input amssym.def
\input amssym.tex

\hsize=12cm
\vsize=18cm
\hoffset=2cm
\voffset=3cm

\footline={\hss{\vbox to 2cm{\vfil\hbox{\rm\folio}}}\hss}

\font\title=cmbx10 at 20pt
\font\head=cmbxti10 at 14pt

\font\mib=cmmib10

\font\tenmib=cmmib10    
\font\sevenmib=cmmib7  
\font\fivemib=cmmib5    
\newfam\mibfam
\def\mib{\fam\mibfam\tenmib}
\textfont\mibfam=\tenmib \scriptfont\mibfam=\sevenmib
\scriptscriptfont\mibfam=\fivemib

\font\teneusm=eusm10    
\font\seveneusm=eusm7  
\font\fiveeusm=eusm5    
\newfam\eusmfam
\def\eusm{\fam\eusmfam\teneusm}
\textfont\eusmfam=\teneusm \scriptfont\eusmfam=\seveneusm
\scriptscriptfont\eusmfam=\fiveeusm

\def\Re{{\rm Re}\,}
\def\Im{{\rm Im}\,}
\def\sgn{{\rm sgn}}
\def\txt#1{{\textstyle{#1}}}
\def\scr#1{{\scriptstyle{#1}}}
\def\f#1{{\goth{#1}}}
\def\r#1{{\rm #1}}
\def\e#1{{\eusm #1}}
\def\varGamma{{\mit\Gamma}}
\def\B#1{{\Bbb #1}}

\vskip 0.5cm
\centerline{\head Spectral Theory of the Riemann Zeta-Function}
\bigskip
\centerline{\head  Chapter 6: Appendix}
\vskip 0.7cm
\hrule
\vskip 0.7cm
\noindent
Leaving the hyperbolic plane $\eusm H$ behind, we now stay on the group
${\rm G}=\r{PSL}(2,{\Bbb R})$. The aim of the present chapter 
is to review what has been related above by changing our vantage point to the
theory of $\varGamma$-automorphic representations of 
$\r{G}$, with $\varGamma=\r{PSL}(2,{\Bbb Z})$. We shall gain, in particular,
a geometrical understanding of the sum formulas involving Kloosterman sums
as well as the explicit formula for the fourth  moment of the 
Riemann zeta-function. We shall obtain also a unified approach to
the mean values of individual automorphic $L$-functions, a subject which
is naturally an extension of the fourth moment of the zeta-function but
does not admit any analogous treatment, thus requiring a genuinely new
method.
\par
We shall develop a relatively self-contained treatment of the spectral theory
of the space $L^2(\varGamma\backslash\r{G})$; in fact, 
this chapter could be read as a particular episode in the theory of unitary
representations of Lie groups. Our reasoning is mostly explicit.
\vskip 1cm
\centerline{\bf 6.1 The group}
\medskip
\medskip
\noindent
To begin with, we introduce a coordinate system into 
$\r{G}$ by means of the Iwasawa decomposition
$$
\r{G}=\r{NAK},\eqno(6.1.1)
$$
where
$$
\eqalignno{
\r{N}&=\left\{\r{n}[x]=\left[\matrix{1&x\cr&1}\right]:
\,x\in{\Bbb R}\right\},\cr
\r{A}&=\left\{\r{a}[y]
=\left[\matrix{\sqrt{y}&\cr&1/\sqrt{y}}\right]:
\,y>0\right\},\cr
\r{K}&=\left\{\r{k}[\theta]=\left[\matrix{\phantom{-}\cos\theta&
\sin\theta\cr-\sin\theta&\cos\theta}\right]:\,\theta\in{\Bbb
R}/\pi{\Bbb Z}\right\}. &(6.1.2)
}
$$
In fact we have
$\r{G}\ni\left[{a\atop c}{b\atop d}\right]=\r{n}[x]\r{a}[y]\r{k}[\theta]$,
with
$$
x={ac+bd\over c^2+d^2},\quad y=(c^2+d^2)^{-1},\quad
\exp(i\theta)={d-ic\over|d-ic|}.
\eqno(6.1.3)
$$
The first two give also that if $\r{n}[x]\r{a}[y]\r{k}[\theta]
=\r{n}[x_1]\r{a}[y_1]\r{k}[\theta_1]$, then $x=x_1$, $y=y_1$, and thus
$\r{k}[\theta]=\r{k}[\theta_1]$ or $\theta\equiv\theta_1\bmod \pi{\Bbb Z}$.
Hence $(6.1.1)$ is indeed a coordinate system on $\r{G}$. We shall always read
the notation $\r{G}\ni\r{g}=\r{n}\r{a}\r{k}$ 
or $\r{n}[x]\r{a}[y]\r{k}[\theta]$ in  this context.
\par
Let $x+iy\in{\eusm H}$ correspond  to the coset $\r{n}[x]\r{a}[y]\r{K}
\in\r{G}/\r{K}$. For $\r{g}=\left[{a\atop c}{b\atop
d}\right]\in\r{G}$,  we have, by $(6.1.3)$,
$\r{g}\r{n}[x]\r{a}[y]=\r{n}[x_1]\r{a}[y_1]\r{k}[\vartheta]$, with 
$$
\eqalign{
x_1={acy^2+(ax+b)(cx+d)\over (cy)^2+(cx+d)^2},\, y_1&=
{y\over (cy)^2+(cx+d)^2},\cr
&\hskip-2.5cm\exp(i\vartheta)=\Big(
{\jmath (\r{g},x+iy)\over|\jmath(\r{g},x+iy)|}\Big)^{-1}.
}\eqno(6.1.4)
$$
The notation $\jmath (\r{g},x+iy)$ is introduced in Section 2.2; here
$\r{g}$ is regarded as an element in ${\Bbb T}(\eusm H)$, that is,
$\r{g}(x+iy)= x_1+iy_1$. Namely, $\r{g}\r{n}[x]\r{a}[y]\r{K}
=\r{n}[x_1]\r{a}[y_1]\r{K}$, and we have
an exact correspondence
between the elements of $\eusm H$ and $\r{G}/\r{K}$ and the one between the
actions of elements of
$\r{G}$ upon them. In this way we identify the pair 
$({\eusm H}, {\Bbb T}({\eusm H}))$ with the pair $(\r{G}/\r{K},\r{G})$.
Note that $(6.1.4)$ is the result of the {\it left\/} multiplication
or translation by $\r{g}$.
\medskip
We then turn to the differentiable structure on $\r{G}$. The coordinate
system $(6.1.1)$ suggests that we should work with 
operators $\partial_x$, $\partial_y$, and $\partial_\theta$. However, 
since the harmonic analysis on $\r{G}$ should contain that on $\eusm H$
which is based on the invariance of the hyperbolic Laplacian $\Delta$,
we need a differentiable structure on $\r{G}$ which commutes with
the left translations by elements of $\r{G}$, 
because the identification of $({\eusm H}, {\Bbb T}({\eusm H}))$ and
$(\r{G}/\r{K},\r{G})$ is built upon those translations. 
A simple and natural way to
realize such a construction is to define the procedure of differentiation on
$\r{G}$ in a manner independent of left translations;  on the other hand the
procedure needs to be the result of group actions 
on $\r{G}$, that is, those done without leaving
$\r{G}$ as the analogy to $\Bbb R$ dictates. 
To satisfy both, it is logical to
utilize {\it right\/} translations. 
\par
To be precise, let us put
$$
{\bf X}_1=\left(\matrix{&1\cr&}\right),\quad
{\bf X}_2=\left(\matrix{1&\cr&-1}\right),\quad
{\bf X}_3=\left(\matrix{&1\cr-1&}\right),\eqno(6.1.5)
$$
and observe that
$$
\eqalignno{
\r{N}&=\left\{\exp({\bf X}_1t):\, t\in{\Bbb R}\right\},\quad
\r{A}=\left\{\exp({\bf X}_2t):\, t\in{\Bbb R}\right\},\cr
\r{K}&=\left\{\exp({\bf X}_3t):\, t\in{\Bbb R}/\pi{\Bbb Z}\right\}.&(6.1.6)
}
$$
Namely $(6.1.1)$ means that these three one-parameter
subgroups generate $\r{G}$ and at each of its element the three {\it curves\/}
pass through. Hence we utilize the differentiation to
the {\it directions\/} ${\bf X}_j$, and define, for any smooth function $f$ on
$\r{G}$,
$$
({\bf x}_jf)(\r{g})=\Big[{d\over dt}\Big]_{t=0}f(\r{g}\exp({\bf X}_jt)).
\eqno(6.1.7)
$$
It is obvious that the differential operators ${\bf x}_j$
are {\it left\/} invariant or commute with any left translation, for they are
defined via {\it right\/} translations. One may use $(6.1.3)$ to see that
at $\r{g}=\r{n}[x]\r{a}[y]\r{k}[\theta]$
$$
\eqalignno{
{\bf x}_1&=y(\cos2\theta)\partial_x+y(\sin2\theta)\partial_y
+(\sin\theta)^2\partial_\theta,\cr
{\bf x}_2&=-2y(\sin2\theta)\partial_x+2y(\cos2\theta)\partial_y
+(\sin2\theta)\partial_\theta,\cr
{\bf x}_3&=\partial_\theta.&(6.1.8)
}
$$
\par
The set  $\{{\bf x}_1,{\bf x}_2, {\bf x}_3\}$ generates 
a non-commutative algebra ${\cal U}$ over $\Bbb C$ with respect to the 
operator multiplication, which is the
set of all left invariant differential operators on $\r{G}$. Its linear subspace
$\f{g}$ generated by  $\{{\bf x}_1,{\bf x}_2, {\bf x}_3\}$ is a Lie
algebra, for we have the commutator relations
$$
[{\bf x}_1,{\bf x}_2]=-2{\bf x}_1,\quad 
[{\bf x}_1,{\bf x}_3]=-{\bf x}_2,
\quad[{\bf x}_2,{\bf x}_3]=4{\bf x}_1-2{\bf x}_3,\eqno(6.1.9)
$$ 
with $[{\bf x}_i,{\bf x}_j]={\bf x}_i{\bf x}_j-{\bf x}_j{\bf x}_i$,
which can be confirmed via $(6.1.8)$; and the Jacobi
identity holds naturally. The space $\f{g}$ is isomorphic to
the Lie algebra generated by the elements $(6.1.5)$, 
for the relations in $(6.1.9)$ hold with ${\bf X}_j$ 
under an obvious correspondence,
as can be verified easily.  
\par
What is important for our purpose is to fix the center 
of ${\cal U}$, and there is a general method to compute an element
in the center. Thus, let the map $\rm{ad}\,{\bf x}$ acting on
$\f{g}$ be defined by $(\rm{ad}\,{\bf x})({\bf y})= [{\bf x},{\bf y}]$. Then 
$$
\hbox{\rm trace of $(\rm{ad}\,{\bf x})\!\cdot\!(\rm{ad}\,{\bf y})$}
$$
is called the Killing form on $\f{g}\times{\f g}$. Let $(k_{ij})$ be
the inverse matrix of the one attached to the form,
with respect to the basis $\{{\bf x}_1,{\bf x}_2, {\bf x}_3\}$. Then the element
$\sum k_{ij}{\bf x}_i {\bf x}_j$ is in the center of $\cal U$.
The relations $(6.1.9)$ give the matrix for the form, and we find
after some rudimentary computation that
$$
\Omega=-{\bf x}_1^2-{1\over4}{\bf x}_2^2+{1\over2}{\bf x}_1{\bf x}_3
+{1\over2}{\bf x}_3{\bf x}_1\eqno(6.1.10)
$$
should be in the center, which one may, however, 
verify directly by using $(6.1.9)$. This is the {\it Casimir\/} operator on
$\r{G}$. 
\par
For the sake of a later purpose we put
$$
\eqalignno{
{\bf w}\;&={\bf x}_3=\partial_\theta,\cr
{\bf e}^+&=2i{\bf x}_1+{\bf x}_2-i{\bf x}_3=e^{2i\theta}(2iy\partial_x
+2y\partial_y-i\partial_\theta),\cr
{\bf e}^-&=-2i{\bf x}_1+{\bf x}_2+i{\bf x}_3=e^{-2i\theta}(-2iy\partial_x
+2y\partial_y+i\partial_\theta).&(6.1.11)
}
$$
The ${\bf e}^\pm$ are termed the Maass operators,
under our normalization. We have the relations
$$
\eqalignno{
[{\bf w},{\bf e}^+]=2i&{\bf e}^+,\quad [{\bf w},{\bf e}^-]=-2i{\bf e}^-,
\quad [{\bf e}^+,{\bf e}^-]=-4i{\bf w},&(6.1.12)\cr
&\Omega=-{1\over4}{\bf e}^+{\bf e}^-
+{1\over4}{\bf w}^2-{1\over2}i{\bf w}.&(6.1.13)
}
$$
From $(6.1.11)$ and $(6.1.13)$ we have also
$$
\Omega=-y^2((\partial_x)^2+(\partial_y)^2)+y\partial_x\partial_\theta.
\eqno(6.1.14)
$$
If a smooth function $f$ on $\r{G}$ is $\r{K}$-trivial, that is,
$f(\r{n}\r{a}\r{k})= f(\r{n}\r{a})$, then we have $\Omega f=\Delta f$.
Since $f$ can be regarded as a function on $\eusm H$, this relation characterizes
the non-Euclidean Laplacian as a restriction of the Casimir operator.
\par
In passing, we remark that $\Omega$ commutes not only with 
left translations but also with right translations. 
In fact, we have, for any smooth $f$,
$$
\int_0^a ({\bf x}_j\Omega f)(\r{g}\exp({\bf X}_jt))dt
=\Omega\int_0^a ({\bf x}_jf)(\r{g}\exp({\bf X}_jt))dt.\eqno(6.1.15)
$$
By the
definition $(6.1.7)$, the left side is equal to 
$(\Omega f)(\r{g}\exp({\bf X}_ja))-(\Omega f)(\r{g})$, while the right side to
$\Omega(f(\r{g}\exp({\bf X}_ja))-f(\r{g}))$. That is,
$(\Omega f)(\r{g}\exp({\bf X}_ja))=\Omega(f(\r{g}\exp({\bf X}_ja)))$,
which gives the assertion. 
\vskip1cm
\centerline{\bf 6.2 Spectral resolution of the Casimir operator}
\medskip
\medskip
\noindent
We now turn to a spectral resolution of the Casimir operator.
First of all we need to fix a measure on $\r{G}$ which generalizes the
non-Euclidean  area element $d\mu$, and is invariant against the left translation.
The latter is naturally required, because of the aforementioned relation between
$({\eusm H},{\Bbb T}(\eusm H))$ and $(\r{G/K}, \r{G})$. We put, in an a priori
manner,
$$
\eqalignno{
d\r{n}&=dx,\,d\r{a}=dy/y,\, d\r{k}=d\theta/\pi,\cr
d\r{g}&=d\r{n}d\r{a}d\r{k}/y= dxdyd\theta/y^2,&(6.2.1)
}
$$
where $\r{g}=\r{n}[x]\r{a}[y]\r{k}[\theta]$, and $dx,\,dy,\,d\theta$ are
ordinary Lebesgue measures. It is immediate that $d\r{n},\,d\r{a},\,
d\r{k}$ are invariant measures on the groups $\r{N},\,\r{A},\,\r{K}$.
As to the invariance of $d\r{g}$, we observe that what is essential,
in the present context, about the change of 
variable $\r{g}\mapsto\r{h}\r{g}$ with a
fixed $\r{h}$ is the nature of the coset map 
$\r{g}\r{K}\mapsto \r{hg}\r{K}$, for
$\r{K}$ is abelian and the measure 
$d\r{k}$ is not affected by the map. 
Thus by the invariance of $d\mu$ we get that of
$d\r{g}$, via $(6.1.4)$. 
\par
The measure $d\r{g}$ is unimodular, i.e., invariant against the 
right translation as well. To verify this, we write $\r{k}[\theta]\r{h}
=\r{n}[\xi(\theta)]\r{a}[u(\theta)]\r{k}[\vartheta(\theta)]$. Then we have
$\r{g}\r{h}=\r{n}[x+\xi(\theta)y]\r{a}[u(\theta)y]
\r{k}[\vartheta(\theta)]$. Thus the Jacobian of the right
translation by $\r{h}$ is equal
to $u(\theta)\vartheta'(\theta)$. An explicit computation of
$\exp(\vartheta(\theta)i)$ via $(6.1.3)$
gives $u(\theta)\vartheta'(\theta)\equiv 1$, which proves the assertion.
In passing, we note further the invariance of $d\r{g}$ against
$\r{g}\mapsto\r{g}^{-1}$, which is, however, irrelevant to our
subsequent discussion.
\medskip
With this, we consider ({\it left\/}) $\varGamma$-automorphic functions $f$ 
on $\r{G}$; that is, $f(\gamma\r{g})=f(\r{g})$ for any pair 
$(\gamma,\r{g})\in\varGamma\!\times\!\r{G}$. If $f$ is smooth,
then ${\bf u}f$ is also $\varGamma$-automorphic for any ${\bf u}\in{\cal U}$,
since $\varGamma$ acts from the left.
The Hilbert space where we work is
$$
L^2(\varGamma\backslash\r{G})=\Big\{
f: \hbox{left $\varGamma$-automorphic and} \int_{\varGamma\backslash\r{G}}
|f(\r{g})|^2d\r{g}<+\infty\Big\},\eqno(6.2.2)
$$
with the natural inner product
$$
\langle{f_1,f_2}\rangle=\int_{\varGamma\backslash G}f_1(\r{g})
\overline{f_2(\r{g})}d\r{g}.
$$
One may choose ${\eusm F}\times\r{K}$ as
$\varGamma\backslash\r{G}$, with a minor abuse of
notation, where ${\eusm F}$ is specified in $(1.1.3)$.
According to Lemma 1.1 and $(6.1.4)$, any coset
$\r{g}\r{K}$ can be mapped by a $\gamma\in\varGamma$ so that 
$\gamma\r{g}\in{\eusm F}\times\r{K}$, and if two
inner points $x+iy$, $x_1+iy_1$ of $\eusm F$ satisfy
$\gamma\r{n}[x]\r{a}[y]\r{k}[\theta]=\r{n}[x_1]\r{a}[y_1]\r{k}
[\theta_1]$ with a $\gamma\in\varGamma$, then we must have $\gamma=1$, 
whence $x=x_1,\,y=y_1,\,\theta\equiv\theta_1\bmod\pi{\Bbb Z}$. 
In fact, one may choose any measurable domain $D$ on $\r{G}$
such that $\gamma D$ ($\gamma\in\varGamma$) cover $\r{G}$ without overlapping
except for sets of null measure. If $D_1$ is another domain of this property,
then for any integrable $\varGamma$-automorphic function $f$
$$
\eqalignno{
\int_{D_1} f(\r{g})d\r{g}&=\sum_{\gamma\in\varGamma}
\int_{\gamma D\cap D_1}f(\r{g})d\r{g}\cr
&=\sum_{\gamma\in\varGamma}
\int_{D\cap\gamma^{-1}D_1}f(\gamma\r{g})d\r{g}=\int_{D} f(\r{g})d\r{g},
&(6.2.3)
}
$$
where the second line is due to the left invariance of $d\r{g}$
and the automorphy of $f$. Hence one may put $\varGamma\backslash\r{G}$
without specifying which domain is under consideration. 
In addition, we note that for $f$ and $D$ as above it holds that
$$
\int_{D}f(\r{gh})d\r{g}=\int_D f(\r{g})d\r{g},\eqno(6.2.4)
$$
with any $\r{h}\in\r{G}$.
For $f(\r{gh})$ is left $\varGamma$ automorphic, $d\r{g}$ is unimodular,
and $D\r{h}$ can obviously stand for $D_1$ in $(6.2.3)$. 
\par
The last observation has an important consequence: With 
smooth vectors $f_1, f_2\in L^2(\varGamma\backslash\r{G})$, 
we have, for any ${\bf u}\in {\cal U}$,
$$
\langle{{\bf u}f_1,f_2}\rangle=\langle{f_1,{\bf u}^*f_2}\rangle,
\eqno(6.2.5)
$$
where
$$
\eqalign{
{\bf u}&=\sum c_{j_1 j_2\ldots j_k}
{\bf x}_{j_1}{\bf x}_{j_2}\cdots{\bf x}_{j_k},\cr
{\bf u}^*&=\sum (-1)^k\overline{c_{j_1 j_2\ldots j_k}}
{\bf x}_{j_k}{\bf x}_{j_{k-1}}\cdots{\bf x}_{j_1},
}
$$
with ${\bf x}_{j_i}\in\{{\bf x}_1, {\bf x}_2,{\bf x}_3\}$.
In fact we have
$$
\eqalignno{
\langle{{\bf x}_jf_1,f_2}\rangle&=\int_{\varGamma\backslash\r{G}}
\left[{d/dt}\right]_{t=0}
f(\r{g}\exp({\bf X}_jt))\overline{f_2(\r{g})}d\r{g}\cr
&=\left[{d/ dt}\right]_{t=0}\int_{\varGamma\backslash\r{G}}
f(\r{g}\exp({\bf X}_jt))\overline{f_2(\r{g})}d\r{g}\cr
&=\left[{d/dt}\right]_{t=0}\int_{\varGamma\backslash\r{G}}
f(\r{g})\overline{f_2(\r{g}\exp(-{\bf X}_jt))}d\r{g}
=\langle{f_1,-{\bf x}_jf_2}\rangle,
}
$$
where the third line is due to $(6.2.4)$. For instance, we have 
$$
\Omega^*=\Omega.\eqno(6.2.6)
$$
Namely, $\Omega$ is self-adjoint.
\medskip
Now, let $f\in L^2(\varGamma\backslash\r{G})$ be smooth and bounded. 
Since $f(\r{g})$ is of period $\pi$ in $\theta$, we have the Fourier expansion
$$
f(\r{g})=\sum_{p=-\infty}^\infty f_p(\r{g})\exp(2pi\theta),\eqno(6.2.7)
$$
with
$$
f_p(\r{g})={1\over\pi}\int_0^\pi f(\r{g}\r{k}[\xi])\exp(-2pi\xi)d\xi.
$$
We may regard the latter as an integral over $\r{K}$, and have,
by the ordinary Parseval identity,
$$
\Vert f(\r{g})\Vert^2_{\r{K}}
=\sum_{p=-\infty}^\infty\Vert f_p(\r{g})\Vert^2_{\r{K}},\eqno(6.2.8)
$$
with the natural norm, which implies the orthonormal decomposition
$$
L^2(\varGamma\backslash\r{G})=\bigoplus_{p=-\infty}^\infty 
L^2_p(\varGamma\backslash\r{G}).\eqno(6.2.9)
$$
Here 
$$
L^2_p(\varGamma\backslash\r{G})=
\left\{f\in L^2(\varGamma\backslash\r{G}): 
f(\r{g}\r{k}[\xi])=f(\r{g})\exp(2pi\xi)\right\}.\eqno(6.2.10)
$$
We then put $g(x+iy)=f(\r{n}[x]\r{a}[y])$ for any $f$
in 
$L^2_p(\varGamma\backslash\r{G})$, and let $\gamma\in\varGamma$ be
such that $\gamma(x+iy)=x_1+iy_1$. We have, by $(6.1.4)$,
$$
\eqalignno{
g(\gamma(x+iy))&=f(\r{n}[x_1]\r{a}[y_1])
=f((\gamma\r{n}[x]\r{a}[y])\cdot\r{k}[-\vartheta])\cr
&=f(\gamma\r{n}[x]\r{a}[y])\exp(-2pi\vartheta)=f(\r{n}[x]\r{a}[y])
\exp(-2pi\vartheta),
}
$$
and thus, with $z=x+iy$,
$$
g(\gamma(z))=g(z)\Big({\jmath(\gamma,z)\over
|\jmath(\gamma,z)|}\Big)^{2p};\eqno(6.2.11)
$$
that is, $g$ is of weight $2p$.
Namely, we  have
$$
L^2(\varGamma\backslash\r{G})=\bigoplus_{p=-\infty}^\infty 
L^2_p(\varGamma\backslash\r{G}/\r{K}),\eqno(6.2.12)
$$
with
$$
\eqalignno{
&L^2_p(\varGamma\backslash\r{G}/\r{K})=
L^2_p(\varGamma\backslash{\eusm H})\cr
&=\Big\{\hbox{$g$ satisfying $(6.2.11)$
and $\displaystyle\int_{\eusm F}|g(z)|^2d\mu(z)
<+\infty$}\Big\}.&(6.2.13)
}
$$
On the other hand,
applying $\Omega$ to both sides of $(6.2.7)$, we have, by $(6.1.14)$,
$$
(\Omega f)(\r{g})=\sum_{p=-\infty}^\infty 
(\Delta_p f_p)(\r{n}[x]\r{a}[y])\exp(2pi\theta),
\eqno(6.2.14)
$$
with $\Delta_p$ as in $(3.2.30)$,
where the smoothness of $f_p$ comes from that of $f$.
This means that the problem of the spectral resolution of $\Omega$ over
$L^2(\varGamma\backslash\r{G})$ is replaced by that of
$\Delta_p$ over $L^2_p(\varGamma\backslash{\eusm H})$. We 
should, however, take a caution.  For, given 
$\Delta_p f_p=\lambda f_p$ with a
certain constant $\lambda$ and with a $p$, 
we are unable to assert immediately that 
$f$ itself is an eigenfunction of $\Omega$ with the eigenvalue $\lambda$,
although the converse is trivial.
We shall show in the sequel that this is in fact the case.
Namely, we shall see that the eigenvalues and eigenvectors of 
$\Delta_p$ over $L^2_p(\varGamma\backslash\r{G}/\r{K})$ with varying $p$
are closely related to each other. Behind this mechanism is the existence of
the Maass operators introduced at $(6.1.11)$.
\medskip
Now, let $f\in L_p^2(\varGamma\backslash\r{G})$ be a $C^2$-class
function such that 
$\Omega f=(\kappa^2+{1\over4})f$; note that $\kappa^2+{1\over4}\in{\Bbb R}$,
because of $(6.2.6)$; here $f$ is such that $f_{p'}\equiv0$ for
$p\ne p'$
in the above notation. We are going to apply differential
operators to $f$, and it is expedient to have an extension of Lemma 1.4. 
Thus, let us put $g(x+iy)=f(\r{n}[x]\r{a}[y])$ as above; we have $\Delta_pg=
(\kappa^2+{1\over4})g$. Since $g$ is of $C^2$-class and of period $1$ 
with respect to $x$, we have the Fourier expansion
$$
g(x+iy)=\sum_{n=-\infty}^\infty b(n,y)e(nx).
$$
Applying $\Delta_p$ to both sides, we see that $b(n,y)$
satisfies the differential equation
$$
-y^2b''(n,y)+((2\pi ny)^2-4\pi npy- \kappa^2-\txt{1\over4})b(n,y)=0.
\eqno(6.2.15)
$$  
If $n\ne0$, then we have, by Lemma 3.8,
$$
b(n,y)=\rho(n)W_{\delta p,i\kappa}(4\pi|n|y),
\quad \delta=\sgn(n),\eqno(6.2.16)
$$
with a certain constant $\rho(n)$. We have used that 
$\Vert{g}\Vert=\Vert{f}\Vert<\infty$, which gives also 
$\rho(n)\ll e^{\varepsilon|n|}$ with any small $\varepsilon>0$,
in view of $(3.2.33)$. 
\par
We shall show that $b(0,y)\equiv0$. The argument
is similar to that in the corresponding part of the proof of Lemma
1.4, though it is a little bit more involved. We may restrict ourselves to the
situation with $i\kappa<0$, and $b(0,y)=\rho(0)y^{{1\over2}+i\kappa}$.
We are going to show that $\rho(0)=0$. To this end, we consider the identity
$$
0=\int_{\e{F}_Y}\left\{\Delta_pg(z)\overline{E_p(z,s)}
-g(z)\overline{\Delta_pE_p(z,s)}\right\}d\mu(z),
$$
where $\e{F}_Y$ is as in $(1.1.30)$, $E_p$ defined by $(3.2.24)$, and
$s={1\over2}-i\kappa$. Note that since we may assume that $p\ne0$, the expansion
$(3.2.27)$ implies that this value of $E_p$ is finite. Integration by parts gives
$$
\eqalign{
0=\int_{\partial\e{F}_Y}\Big({\partial g\over\partial\tilde{n}}(z)
\overline{E_p(z,s)}
-g(z){\partial\over\partial\tilde{n}}\overline{E_p(z,s)}
+2ipg(z)\overline{E_p(z,s)}\Big){|dz|\over y},
}
$$
where $\partial/\partial\tilde{n}$ is the non-Euclidean outer-normal 
differentiation introduced in Section 1.1. Let $\gamma\in\varGamma$. 
Since $\gamma$ commutes with $\partial/\partial\tilde{n}$, we have,
by $(3.2.26)$ and $(6.2.11)$,
$$
\eqalign{
{\partial g\over\partial\tilde{n}}&(\gamma(z))\overline{E_p(\gamma(z),s)}=
{\partial\gamma g\over\partial\tilde{n}}(z)\overline{E_p(\gamma(z),s)}\cr
&={\partial g\over\partial\tilde{n}}(z)\overline{E_p(z,s)}
+g(z)\overline{E_p(z,s)}\Big({\jmath(\gamma,z)\over|\jmath(\gamma,z)|}\Big)^{-2p}
{\partial\over\partial\tilde{n}}
\Big({\jmath(\gamma,z)\over|\jmath(\gamma,z)|}\Big)^{2p}.
}
$$
Putting $\jmath(\gamma,z)=cz+d$, one may compute the last derivative
explicitly. We have
$$
\eqalign{
\Big({\jmath(\gamma,z)\over|\jmath(\gamma,z)|}\Big)^{-2p}
{\partial\over\partial\tilde{n}}
\Big({\jmath(\gamma,z)\over|\jmath(\gamma,z)|}\Big)^{2p}&=
-2piy{d\over|dz|}\log|\jmath(\gamma,z)|\cr
&=piy{d\over|dz|}\log\Big({\Im\gamma(z)\over\Im z}\Big).
}
$$
Collecting these, we find that
$$
\eqalign{
0=\int_{-{1\over2}}^{1\over2}
\Big({\partial g\over\partial y}(z)
\overline{E_p(z,s)}
-g(z){\partial\over\partial y}\overline{E_p(z,s)}
+2ipy^{-1}g(z)\overline{E_p(z,s)}\Big)_{y=Y}dx,
}
$$
which gives $\rho(0)=0$. 
\medskip
That is, we have, for any $f\in L^2_p(\varGamma\backslash\r{G})$ such that 
$\Omega f=(\kappa^2+{1\over4})f$,
$$
f(\r{g})=e^{2pi\theta}
\sum_{\scr{n=-\infty}\atop\scr{n\ne0}}^\infty\rho(n)
W_{\delta p,i\kappa}(4\pi|n|y)e(nx)
$$
With this, let us consider $f^{(-1)}={\bf e}^-f$.   
We have, by $(6.1.11)$,
${\bf w}f^{(-1)}={\bf e}^-{\bf w}f-2i{\bf e}^-f={2i(p-1)}f^{(-1)}$. 
Moreover, by $(6.1.13)$ and $(6.2.5)$
$$
\eqalignno{
\Vert{f^{(-1)}}\Vert^2&=\langle{{\bf e}^-f,{\bf e}^{-}f}\rangle
=\langle{f,-{\bf e}^+{\bf e}^{-}f}\rangle\cr
&=\langle{f,(4\Omega-{\bf w}^2+2i{\bf w})f}\rangle=
4(\kappa^2+(p-\txt{1\over2})^2)\Vert{f}\Vert^2,&(6.2.17)
}
$$ 
whence
$f^{(-1)}\in L^2_{p-1}(\varGamma\backslash G)$; moreover,
$\Omega f^{(-1)}=(\kappa^2+{1\over4})f^{(-1)}$, for $\Omega$ 
is in the center of $\cal U$. Analogously, with $f^{(+1)}={\bf e}^+f$,
we have $f^{(+1)}\in L^2_{p+1}(\varGamma\backslash G)$, 
$\Omega f^{(+1)}=(\kappa^2+{1\over4})f^{(+1)}$ as well as
$$
\Vert{f^{(+1)}}\Vert^2=4(\kappa^2+(p+\txt{1\over2})^2)\Vert{f}\Vert^2.
\eqno(6.2.18)
$$
We may repeat the procedure like ascending and descending 
aerial strata. There are three possible cases:
\medskip
\item{(1)} $\kappa\ge0$, 
\item{(2)} $\Im\kappa>0$ but $\kappa\ne i(q-{1\over2})$ for any
integer $q$,
\item{(3)} $\kappa=i(q-{1\over2})$, with an integer $q>0$.
\medskip
\noindent
In both the cases (1) and (2) with $\pm p>0$ we have that $f^{(\mp p)}
\not\equiv0$ is $\r{K}$-trivial, i.e., in $L^2(\varGamma\backslash{\eusm H})$. 
Namely, there exists a real analytic cusp form $\psi$ such that
$\Delta\psi=(\kappa^2+{1\over4})\psi$ and 
$f^{(\mp p)}=c\psi$ with a constant $c$. Then Lemma 1.4 implies that
the case (2) is impossible under our assumption that $\varGamma
=\r{PSL}(2,{\Bbb Z})$. Also, with (1) we have in fact $\kappa>3.815$; and
the procedure can be reversed. We may express 
this fact that $f$ can be reached
by either ascending or descending from a real analytic cusp form $\psi$, i.e., 
$({\bf e}^\pm)^p\psi=f$. 
\par
On the other hand, in the case (3) with $p>0$ the descent
terminates, for we have $f^{(-p+q-1)}\equiv0$ as $(6.2.17)$
implies. The Fourier coefficients of $f^{(-p+q)}$ satisfy $(6.2.15)$
with $p=q$ and $\kappa^2+{1\over4}=-q(q-1)$. We may use the second identity
in $(3.2.35)$. On noting that $W_{-q,q-{1\over2}}(y)$ 
satisfies the same
differential equation as that for $W_{q,q-{1\over2}}(-y)$,
we have
$$
f^{(-p+q)}(\r{g})
=y^qe^{2iq\theta}\sum_{\scr{n=-\infty}\atop\scr{n\ne0}}^\infty
\rho(n)\exp(-2\pi|n|y+2\pi inx).
$$
The last sum, denoted by
$h(z)$, $z=x+iy$, should converge absolutely for any $z\in{\eusm F}$.
The equation $\Omega f^{(-p+q)}=-q(q-1)f^{(-p+q)}$ implies that
$$
\Delta h(z)={8\pi^2q  y}
\sum_{n<0}|n|\rho(n)\exp(-2\pi|n|y+2\pi inx).
$$
This gives $\rho(n)=0$ for all $n<0$, since $\Delta h\equiv0$. 
Hence we have found that
$$
f^{(-p+q)}(\r{g})
=y^qe^{2iq\theta}\sum_{n=1}^\infty\rho(n)e(nz).
$$
By the same way as the derivation of $(6.2.11)$ we have
$h(\gamma(z))=(\jmath(\gamma,z))^{2q}h(z)$ for any $\gamma\in\varGamma$;
that is, $h(z)$ is a holomorphic cusp form of weight $2q$ with respect
to $\varGamma$. The case (2) with $p<0$ is analogous, and the counterpart
of $h$ turns out to be anti-holomorphic, i.e.,
$$
f^{(-p-q)}(\r{g})
=y^qe^{-2iq\theta}\sum_{n=1}^\infty \rho(n)e(-n\overline{z}),
$$
in which the complex conjugate of the sum is a holomorphic cusp
form of weight $2q$.
\medskip
What remains then is to make precise the
contribution of the continuous spectrum. This can also be dealt with in a fashion
similar to the above; that is, the action of the Maass operators is 
again the key. Thus, from what we have seen in 
the above, we expect that except for
those vectors originating from holomorphic cusp forms
the space $L^2_p(\varGamma\backslash\r{G})$, $p\ge0$, should be spanned by the 
$({\bf e}^+)^p$-images of smooth vectors of
$L^2_0(\varGamma\backslash\r{G})= L^2(\varGamma\backslash{\eusm H})$.
\par
In order to verify this proposition, we first re-define 
the Eisenstein series $E_p$
introduced at $(3.2.24)$. We put, for any $p\in\Bbb Z$,
$$
\phi_p(\r{g},\nu)=y^{\nu+{1\over2}}e^{2pi\theta},\eqno(6.2.19)
$$
and
$$
E_p(\r{g},\nu)=\sum_{\gamma\in\varGamma_\infty\backslash\varGamma}
\phi_p(\gamma\r{g},\nu),\quad \Re\nu>{1\over2}.\eqno(6.2.20)
$$
Via $(6.1.4)$, we have that
$$
E_p(\r{g},\nu)=E_p(z,\nu+\txt{1\over2})e^{2pi\theta},\quad z=x+iy,
\eqno(6.2.21)
$$
where the $E_p$ on the right side stands for $(3.2.24)$; there should
not be any notational confusion. On noting that
$$
{\bf e}^\pm\phi_p(\r{g})=(2\nu+1\pm2p)\phi_{p\pm1}(\r{g}),\eqno(6.2.22)
$$
we have, for $p\ge0$,
$$
({\bf e^\pm})^p E_0(\r{g},\nu)=\prod_{\ell=0}^{p-1}
(2\nu+1\pm2\ell)\cdot
E_{\pm p}(\r{g},\nu).\eqno(6.2.23)
$$
\par
Then, let us consider an $f$ in $L^2_p(\varGamma\backslash\r{G})$
with $p>0$, each partial derivative of which is of fast decay; c.f.,
$(1.1.29)$. Since $({\bf e}^{-})^pf\in L^2(\varGamma\backslash{\eusm H})$,
we have, by Theorem 1.1 with a minor rearrangement, 
$$
\eqalignno{
({\bf e}^{-})^pf(\r{g})
&=\sum_{j=0}^\infty\langle{({\bf e}^{-})^pf,
\psi_j}\rangle\psi_j(\r{g})\cr
&+{1\over2\pi}\int_0^\infty E_0(\r{g},it)
{\eusm E}(t,({\bf e}^{-})^pf)dt,&(6.2.24)
}
$$
By $(6.2.23)$ we rewrite this as
$$
\eqalignno{
({\bf e}^{-})^pf(\r{g})
&=(-1)^p\sum_{j=0}^\infty\langle{f,
({\bf e}^{+})^p\psi_j}\rangle\psi_j(\r{g})\cr
&+{(-1)^p\over2\pi}\int_0^\infty
E_0(\r{g},it) {\eusm E}_p(t,f)dt,&(6.2.25)
}
$$
with
$$
{\eusm E}_p(t,f)=2^p{\Gamma({1\over2}-it+p)\over
\Gamma({1\over2}-it)}\int_{\varGamma\backslash\r{G}}f(\r{g})\overline{
E_p(\r{g},it)}d\r{g}.\eqno(6.2.26)
$$
We then observe that $(6.1.12)$, $(6.1.13)$, and $(6.2.23)$ give
$$
\eqalignno{
({\bf e}^-)^p({\bf e}^+)^p\psi_j(\r{g})
&=(-4)^p{|\Gamma({1\over2}+i\kappa_j+p)|^2
\over|\Gamma({1\over2}+i\kappa_j)|^2}\cdot\psi_j(\r{g}),\cr
({\bf e}^-)^p({\bf e}^+)^pE_0(\r{g},it)&=
(-2)^p{\Gamma({1\over2}+it+p)
\over\Gamma({1\over2}+it)}\cdot ({\bf e}^-)^pE_p(\r{g},it). &(6.2.27)
}
$$
A combination of $(6.2.25)$--$(6.2.27)$ yields that
$({\bf e}^-)^p f^*(\r{g})=0$,
with
$$
\eqalignno{
f^*(\r{g})=f(\r{g})
&-2^{-2p}\sum_{j=0}^\infty\langle{f,({\bf e}^+)^p\psi_j}\rangle
{|\Gamma({1\over2}+i\kappa_j)|^2
\over|\Gamma({1\over2}+i\kappa_j+p)|^2}\cdot({\bf e}^+)^p\psi_j(\r{g})\cr
&-{1\over2\pi}\int_0^\infty E_p(\r{g},it){\cal E}_p(t,f)dt.&(6.2.28)
}
$$
We assert that
$$
f^*=\sum_{\ell=0}^p ({\bf e}^+)^{p-\ell}\varphi_\ell,\eqno(6.2.29)
$$
where $y^{-\ell}\varphi_\ell(\r{n}[x]\r{a}[y])$
is a holomorphic cusp form of weight $2\ell$; note that we have
actually $\ell\ge6$, since there exist no holomorphic cusp forms
of weight less than $12$ over $\varGamma$. 
We prove $(6.2.29)$ by induction with respect to $p$. Thus, let
$f_1\in L^2_{p+1}(\varGamma\backslash\r{G})$ be such that 
$({\bf e}^-)^{p+1}f_1\equiv0$. By the inductive assumption, we have
${\bf e}^-f_1=\sum_{\ell=0}^p ({\bf e}^+)^{p-\ell}\varphi_{1,\ell}$,
where the specification of the right side is as that of $(6.2.29)$. Applying
${\bf e}^+$ to both sides, we have, by $(6.1.13)$,
$$
\eqalignno{
(\Omega+p(p+1))f_1&=-{1\over4}
\sum_{\ell=0}^p ({\bf e}^+)^{p+1-\ell}\varphi_{1,\ell}\cr
&=-{1\over4}\sum_{\ell=0}^p ({\bf e}^+)^{p+1-\ell}
{(\Omega+p(p+1))\over
\ell(1-\ell)+p(p+1)}\varphi_{1,\ell}.
}
$$
That is, $(\Omega+p(p+1))f_2=0$ with
$$
f_2=f_1+{1\over4}\sum_{\ell=0}^p 
{({\bf e}^+)^{p+1-\ell}\varphi_{1,\ell}\over
\ell(1-\ell)+p(p+1)}.
$$
Thus $y^{-p-1}f_2(\r{n}[x]\r{a}[y])$ is a holomorphic cusp form
of weight $2(p+1)$, which ends the proof of $(6.2.29)$.
\medskip
The formula $(6.2.28)$ with $(6.2.29)$ reveals the spectral
structure of the space $L^2_p(\varGamma\backslash\r{G})$,
$p>0$. It is generated by the $({\bf e}^+)^p$-images 
of real analytic cusp forms and integrals of 
Eisenstein series  and by the vectors
of  the type $(6.2.29)$. The case $p<0$ is analogous.
\vskip1cm
\centerline{\bf 6.3 Automorphic representations}
\medskip
\medskip
\noindent
The above discussion essentially completes 
the spectral resolution of $\Omega$ over $L^2(\varGamma\backslash \r{G})$. 
We see horizontal strata in the space, 
all of which whirl with individual rates by
the action of $\r{K}$ from the right; thus, no mixing takes
place.  There are, however, lifts to
climb up and down the strata.  Their vertical ways naturally never cross each
other,  that is, they are invariant against the action of $\r{G}$ from the right,
as is implicitly asserted at $(6.1.15)$. Some penetrate the hyperbolic plane, and 
some start or terminate without touching it. 
We are now about to render this structure in terms of the
$\varGamma$-automorphic representation of the Lie group $\r{G}$.
\medskip
To begin with, we rearrange the somewhat complicated formula $(6.2.32)$ by using
a certain integral transform due to H. Jacquet.
As an orientation, we observe that in view of the Fourier expansion $(3.2.27)$
of $E_p$ on the right side of $(6.2.24)$ it is natural to expect that there
should exist an integral transform connecting $\phi_p(\r{g},\nu)$ with
the Whittaker function $W_{p,\nu}(y)$, which defines the latter in terms
of elements and actions of $\r{G}$, and makes it possible to
understand the basic differential equation $(3.2.32)$ in terms of
the pair $(\r{G},\f{g})$.
\par
The integral transform we are concerned is defined by
$$
{\cal A}^\delta\!f(\r{g})=\int_{-\infty}^\infty
e(-\delta\xi)f(\r{w}\r{n}[\xi]\r{g})d\xi,\eqno(6.3.1)
$$
whenever the integral converges absolutely; here $\delta=\pm1$, and
$\r{w}=\left[{\atop-1}{1\atop}\right]$ is the Weyl element of $\r{G}$. 
A basic property of ${\cal A}^\delta$ is that it commutes
with right translations, and is an inter-twining
operator; that is, we have, for any ${\bf u}\in{\cal U}$,
$$
{\bf u}({\cal A}^\delta\!f)(\r{g})={\cal A}^\delta({\bf u}f)(\r{g}),\eqno(6.3.2)
$$
provided $f$ is smooth. In view of $(6.1.7)$, 
this is immediate with ${\bf u}={\bf x}_j$, and the general case
as well. As we shall see in Section 6.5
below, ${\cal A}^\delta$ is closely related 
to the Fourier expansion of Poincar\'e
series on $\r{G}$, with respect to the left action of $\r{N}$.
\medskip
By $(6.1.3)$ the map $\r{g}\mapsto\r{w}\r{n}[\xi]\r{g}$ is equivalent to
$$
\eqalignno{
&x\mapsto {-x-\xi\over\sqrt{y^2+(x+\xi)^2}},\quad
y\mapsto {y\over y^2+(x+\xi)^2}\cr
&\hskip 1cm e^{2pi\theta}\mapsto e^{2pi\theta}\Big({x+\xi-iy\over
x+\xi+iy}\Big)^p.&(6.3.3)
}
$$
Hence
$$
\eqalignno{
{\cal A}^\delta\!\phi_p(\r{g},\nu)&=\exp(2pi\theta)e(\delta x)
y^{-\nu+{1\over2}}\int_{-\infty}^\infty{e(y\xi)
\over(\xi^2+1)^{\nu+{1\over2}}}\Big({\xi+i\over\xi-i}\Big)^{\delta p}d\xi\cr
&=(-1)^p\pi^{\nu+{1\over2}}\exp(2pi\theta)e(\delta x)
{W_{\delta p,\nu}(4\pi y)\over\Gamma(\delta p+\nu+{1\over2})}.&(6.3.4)
}
$$
The first line is valid for $\Re\nu>0$, and the second line for all
$\nu\in{\Bbb C}$ because of $(3.2.31)$ and Lemma 3.8.
On noting the first identity in $(3.2.35)$,
this implies that we may rewrite $(1.1.43)$
as
$$
\eqalignno{
\psi_j(\r{g})&={\Gamma(\txt{1\over2}+i\kappa_j)\over
2\pi^{{1\over2}+i\kappa_j}}
\sum_{\scr{n=-\infty}\atop
\scr{n\ne0}}^\infty{\rho_j(n)\over\sqrt{|n|}}
{\cal A}^{\sgn(n)}\!\phi_0(\r{a}[|n|]\r{g},i\kappa_j)\cr
&=\sum_{\scr{n=-\infty}\atop
\scr{n\ne0}}^\infty{\varrho_j(n)\over\sqrt{|n|}}
{\cal A}^{\sgn(n)}\!\phi_0(\r{a}[|n|]\r{g},i\kappa_j),&(6.3.5)
}
$$
with the new normalization of the Fourier coefficients:
$$
\varrho_j(n)={\Gamma(\txt{1\over2}+i\kappa_j)\over
2\pi^{{1\over2}+i\kappa_j}}\rho_j(n).\eqno(6.3.6)
$$
By virtue of $(6.3.2)$, we have
$$
({\bf e}^+)^p\psi_j(\r{g})=
\sum_{\scr{n=-\infty}\atop
\scr{n\ne0}}^\infty{\varrho_j(n)\over\sqrt{|n|}}
{\cal A}^{\sgn(n)}\!({\bf e}^+)^p\phi_0(\r{a}[|n|]\r{g},i\kappa_j).
\eqno(6.3.7)
$$
Then we put
$$
\lambda_j^{(p)}(\r{g})={\Gamma({1\over2}+i\kappa_j)
\over 2^p\Gamma(\txt{1\over2}+i\kappa_j+p)}({\bf e}^+)^p\psi_j(\r{g}).
\eqno(6.3.8)
$$
We see readily that
$$
\lambda_j^{(p)}(\r{g})=\sum_{\scr{n=-\infty}\atop
\scr{n\ne0}}^\infty{\varrho_j(n)\over\sqrt{|n|}}
{\cal A}^{\sgn(n)}\!\phi_p(\r{a}[|n|]\r{g},i\kappa_j).\eqno(6.3.9)
$$
We have
$$
\langle{\lambda_j^{(p)},\lambda_l^{(p)}}\rangle=\delta_{jl}.\eqno(6.3.10)
$$
In fact, this is the same as the first formula in $(6.2.27)$. Hence, we
may rewrite $(6.2.28)$ as
$$
f(\r{g})=f^*(\r{g})
+\sum_{j=0}^\infty\langle{f,\lambda_j^{(p)}}\rangle
\lambda_j^{(p)}(\r{g})
+{1\over2\pi}\int_0^\infty E_p(\r{g},it){\cal E}_p(t,f)dt,\eqno(6.3.11)
$$
for any smooth $f\in L_p^2(\varGamma\backslash\r{G})$.
\par
We extend $(6.3.9)$ to those vectors of $L_p^2(\varGamma\backslash\r{G})$ which
come from holomorphic cusp forms. Thus, let ${\cal{C}}_k(\varGamma)$ with
$k\le p$ be the vector space of holomorphic cusp forms of weight $2k$ that is
introduced in Section 2.2. Let $\{\psi_{j,k}\,:\, 1\le{j}\le\vartheta(k)\}$ be
its orthonormal base defined in $(2.2.2)$. We then put
$$
\lambda_{j,k}^{(k)}(\r{g})=y^k\psi_{j,k}(x+iy)\exp(2ki\theta).\eqno(6.3.12)
$$
Also, corresponding to $(6.3.6)$, we introduce the renormalization 
of the Fourier coefficients of $\psi_{j,k}$:
$$
\varrho_{j,k}(n)=(-1)^k{\Gamma(2k)^{1\over2}\over 2^{2k}\pi^{k+{1\over2}}}
\rho_{j,k}(n).\eqno(6.3.13)
$$
By the first line of $(6.3.4)$, we may rewrite $(2.2.3)$ as
$$
\lambda_{j,k}^{(k)}(\r{g})=\pi^{{1\over2}-k}\Gamma(2k)^{1\over2}
\sum_{n=1}^\infty{\varrho_{j,k}(n)\over\sqrt{n}}
{\cal A}^+\!\phi_k(\r{a}[n]\r{g},k-\txt{1\over2}),\eqno(6.3.14)
$$
which is a counterpart of $(6.3.5)$. Further, we put 
$$
\lambda_{j,k}^{(p)}(\r{g})=2^{k-p}
\Big({\Gamma(2k)\over\Gamma(p-k+1)\Gamma(p+k)}\Big)^{1\over2}({\bf
e}^+)^{p-k}\lambda_{j,k}^{(k)}(\r{g}).\eqno(6.3.15)
$$
We have
$$
\lambda_{j,k}^{(p)}(\r{g})=\pi^{{1\over2}-k}
\Big({\Gamma(p+k)\over\Gamma(p-k+1)}\Big)^{1\over2}
\sum_{n=1}^\infty{\varrho_{j,k}(n)\over\sqrt{n}}
{\cal A}^+\!\phi_p(\r{a}[n]\r{g},k-\txt{1\over2}).\eqno(6.3.16)
$$
As an analogue of $(6.3.10)$, we have, for any $p\ge k$,
$$
\langle{\lambda_{j,k}^{(p)},\,\lambda_{l,k}^{(p)}}\rangle=\delta_{j,l}.
\eqno(6.3.17)
$$
Hence $(6.2.33)$ can be expressed as
$$
f^*(\r{g})=\sum_{k=6}^p\sum_{j=1}^{\vartheta(k)}
\langle{f,\lambda_{j,k}^{(p)}}\rangle \lambda_{j,k}^{(p)}(\r{g}).\eqno(6.3.18)
$$
In fact, $f^*$ is obviously in the space spanned by 
$\{\lambda_{j,k}^{(p)}:\,1\le j\le\vartheta(k),\,k\le p\}$, which is
orthogonal to the space spanned by $\{\lambda_j^{(p)}:\, j=1,2,\ldots,\infty\}$
and integrals of the Eisenstein series $E_p$. For instance,
$\langle\lambda_j^{(p)},\,\lambda_{l,k}^{(p)}\rangle=0$ follows from
the identity $\langle\Omega\lambda_j^{(p)},\,\lambda_{l,k}^{(p)}\rangle=
\langle\lambda_j^{(p)},\,\Omega\lambda_{l,k}^{(p)}\rangle$;  
the Eisenstein series is treated analogously.
\medskip
The space $L_p^2(\varGamma\backslash\r{G})$ with $p<0$ has essentially the
same spectral structure, except that anti-holomorphic cusp
forms fill the r\^ole of holomorphic cusp forms. Their orthonormal base
vectors are given by
$$
\lambda_{j,k}^{(p)}(\r{g})=\pi^{{1\over2}-k}
\Big({\Gamma(|p|+k)\over\Gamma(|p|-k+1)}\Big)^{1\over2}
\sum_{n=1}^{\infty}{\varrho_{j,k}(n)\over\sqrt{n}}
{\cal A}^-\!\phi_p(\r{a}[n]\r{g},k-\txt{1\over2}),\eqno(6.3.19)
$$
with $\varrho_{j,k}(n)$ as in $(6.3.16)$. The involution $\r{g=nak}\mapsto
\r{n}^{-1}\r{a}\r{k}^{-1}$ sends anti-holomorphic cusp forms to
holomorphic cusp forms, and vice versa.
\medskip
\medskip
Collecting the above discussion, we obtain
\medskip
\medskip
\noindent
{\bf Theorem 6.1} (The spectral resolution of the Casimir operator)
{\it Let ${}^0\!L^2(\varGamma\backslash\r{G})$ be 
the cuspidal subspace spanned by those vectors 
in $L^2(\varGamma\backslash\r{G})$
whose Fourier expansion with respect to 
the left action of $\r{N}$ have vanishing
constant terms, and 
${}^e\!L^2(\varGamma\backslash\r{G})$ the subspace spanned
by integrals of Eisenstein series
$E_p$ of all even integral weights which are defined by $(6.2.20)$. 
Then we have the orthogonal decomposition
$$
L^2(\varGamma\backslash\r{G})={\Bbb C}\cdot1\oplus
{}^0\!L^2(\varGamma\backslash\r{G})\oplus{}^e\!L^2(\varGamma\backslash\r{G}).
\eqno(6.3.20)
$$
More precisely, we have the orthonormal decomposition
$$
{}^0\!L^2(\varGamma\backslash\r{G})=\bigoplus V,\eqno(6.3.21)
$$
where $V$ runs over all
$$
V_j=\bigoplus_{p=-\infty}^\infty {\Bbb C}\lambda_j^{(p)},\quad
V_{l,k}^\pm=\bigoplus_{\pm p=k}^\infty {\Bbb C}\lambda_{l,k}^{(p)}, 
\eqno(6.3.22)
$$
with the base elements being defined by $(6.3.9)$, $(6.3.16)$, and
$(6.3.19)$, respectively; thus $j$ varies from $1$ to infinity, $k$ from $6$
to infinity, and $1\le l\le\vartheta(k)$, with $\vartheta(k)$ as in $(2.2.2)$.
Also we have the orthonormal decomposition
$$
{}^e\!L^2(\varGamma\backslash\r{G})=\bigoplus_{p=-\infty}^\infty
E^{(p)},\eqno(6.3.23)
$$
where $E^{(p)}$ consists of
$$
{1\over4\pi i}\int_{(0)}E_p(\r{g},\nu)h(\nu)d\nu,\eqno(6.3.24)
$$
with $h$ being ordinary square integrable functions over the imaginary
axis $(0)$.
}
\medskip
\medskip
\noindent
A few points are missing in our discussion so far developed. We have started
with $(6.2.7)$; and each component $f_p$ there has been spectrally expanded
as in $(6.3.11)$ together with $(6.3.18)$, although
there is a slight notational confusion. We have not mentioned the
detail about several convergence issues, for their verifications are 
immediate as far as $f$ on $\r{G}$ is left $\varGamma$-automorphic and such that
${\bf u}f$ with ${\bf u}\in{\cal U}$
of sufficiently high order are all of fast decay; 
and the set of those $f$ is 
dense in $L^2(\varGamma\backslash\r{G})$. It would be expedient to remark
here that the {\it uniform\/} fast decay of projections of $f$ to any of the
subspaces listed in $(6.3.22)$ and $(6.3.23)$ 
could be derived readily from that of ${\bf u}f$; for instance,
$$
\eqalignno{
\langle{f,\lambda_j^{(p)}}\rangle&=
{1\over(\kappa_j^2+\txt{1\over4}+2pi)^a}
\langle{f,(\Omega-{\bf
w})^a \lambda_j^{(p)}}\rangle\cr
&={1\over(\kappa_j^2+\txt{1\over4}+2pi)^a}
\langle{(\Omega+{\bf w})^af, \lambda_j^{(p)}}\rangle\cr
&\ll (\kappa_j^2+|p|)^{-a}\Vert{(\Omega+{\bf w})^af}\Vert,
&(6.3.25)
}
$$
with any fixed integer $a\ge0$, gives what is needed, and the same device
works for all other subspaces. We have neither 
mentioned explicitly the Parseval formula 
which generalizes $(1.1.49)$ to the whole
$L^2(\varGamma\backslash\r{G})$; however, this could readily be inferred
from a combination of $(6.2.7)$, $(6.2.8)$, and $(6.3.11)$, since the last is
in fact a rearrangement of $(6.2.24)$, which is in turn a direct
consequence of $(1.1.47)$.
\medskip
\medskip
We are now at the stage to express the assertion of the last theorem
in the language of automorphic representations. Thus, let us consider
the right translation
$$
\omega(h):\,f(\r{g})\mapsto f(\r{g}\r{h}),\eqno(6.3.26)
$$
with any $f\in L^2(\varGamma\backslash\r{G})$. For each $h\in\r{G}$, $\omega(h)$
is a unitary map of $L^2(\varGamma\backslash\r{G})$ into itself, because of 
the unimodularity of $d\r{g}$; and $\omega$ is a homomorphism. 
This configuration is expressed that
$\omega$ is a $\varGamma$-automorphic unitary representation of the Lie group
$\r{G}$. 
\medskip
If $W$ is a closed subspace of $L^2(\varGamma\backslash\r{G})$ and
$\omega(\r{h})W\subset W$ for all $\r{h}\in\r{G}$, then $W$ is 
called an invariant subspace. 
We shall prove that those $V$ in $(6.3.21)$ are all
invariant subspaces.  Thus, let  
$\tilde{V}$ be the closed subspace generated by $\bigcup_\r{h}\omega(\r{h})V$
with $\r{h}$  varying throughout $\r{G}$. Since $\Omega$ commutes
with any right translation as is shown at $(6.1.15)$, all smooth elements in
$\tilde{V}$ are eigenfunctions of
$\Omega$ with the same eigenvalue. The Fourier coefficients, with respect
to the right action of $\r{K}$, of a particular
eigenfunction come from cusp forms either real analytic or holomorphic or
anti-holomorphic over
${\eusm H}$; hence $\tilde{V}$ splits into a finite number of subspaces among
those listed in $(6.3.22)$ which share the same eigenvalue 
with $V$. This means that $\r{G}$ splits into 
the same number of cosets $\r{H}\r{h}$
with $\r{H}$ being the closed subgroup of
$\r{G}$ composed of all elements that send $V$ to itself.
Let us assume that this number is larger than $1$, and consider the
cosets $\r{H}\exp(t{\bf X}_j)$, $t\in\Bbb R$. With either $j=1$ or $2$, 
there should be at least two different cosets. This is, however,
a contradiction, for $\exp(t{\bf X}_j)$ is of course a continuous 
curve on $\r{G}$. 
\medskip
Those $V$'s are in fact irreducible representations; that is, any
invariant subspace contained in $V$ is either $V$ itself or
$\{0\}$. The proof of this fact requires some preparation which
appears to be an excess for our present purpose. 
In fact, what we need genuinely is not the irreducibility but
a realization of the structure of each $V$ in terms of
an ordinary functional space that we shall develop in the
next section. Nevertheless, we shall see that the latter gives 
a somewhat unconventional proof of the former as well. Thus 
in the closing paragraph of the next section we shall prove
\medskip
\medskip
\noindent
{\bf Theorem 6.2} {\it The identity $(6.3.21)$  gives the decomposition
of the cuspidal subspace into irreducible subspaces with respect to
the unitary representation $(6.3.26)$ of $\r{G}$.
}
\medskip
\medskip
Representations and invariant subspaces are obviously 
inter-changeable concepts. With this convention, we may 
call $V$'s in $(6.3.21)$ as irreducible representations of $\r{G}$ 
occurring in the Hilbert space
$L^2(\varGamma\backslash\r{G})$. Those $V_j$ arising from real analytic
cusp forms belong to {\it the unitary principal series\/} of
irreducible representations, and $V_{l,k}^\pm$ 
to {\it the holomorphic\/}
and {\it the anti-holomorphic discrete series\/}, respectively. In general
there can be additional series of representations coming from
exceptional eigen vectors of the Casimir operator ({\it the complementary
series\/}); in our situation with $\varGamma$, such  
does not occur.
\medskip
Here we introduce a major simplification of notation. We write
$$
\eqalignno{
\lambda_V^{(p)}(\r{g})&=\sum_{\scr{n=-\infty}\atop
\scr{n\ne0}}^\infty{\varrho_V(n)\over\sqrt{|n|}}
{\cal A}^{\sgn(n)}\!\phi_p(\r{a}[|n|]\r{g},\nu_V),&(6.3.27)\cr
\lambda_V^{(p)}(\r{g})&=\pi^{-\nu_V}
\Big({\Gamma(|p|+\nu_V+{1\over2})
\over\Gamma(|p|-\nu_V+{1\over2})}\Big)^{1\over2}
\sum_{n=1}^\infty{\varrho_V(n)\over\sqrt{n}}
{\cal A}^\pm\!\phi_p(\r{a}[n]\r{g},\nu_V),\quad&(6.3.28)\cr
E_p(\r{g},\nu)&=y^{{1\over2}+\nu}\exp(2pi\theta)
+(-1)^p{\Gamma({1\over2}+\nu)^2\varphi_\varGamma({1\over2}+\nu)\over
\Gamma({1\over2}+\nu+p)\Gamma({1\over2}+\nu-p)}\exp(2pi\theta)\cr
&\qquad+\sum_{\scr{n=-\infty}\atop
\scr{n\ne0}}^\infty{|n|^{-\nu}\sigma_{2\nu}(|n|)\over\zeta(1+2\nu)\sqrt{|n|}}
\e{A}^{\sgn(n)}\phi_p(\r{a}[|n|]\r{g},\nu).&(6.3.29)
}
$$
The first is equivalent to $(6.3.9)$, the second to either $(6.3.16)$
or $(6.3.19)$, and the third to $(3.2.27)$. 
This is possible, for the Fourier coefficients
$\varrho_j(n)$  and $\varrho_{j,k}(n)$ do not depend on the weight $2p$ but
only on the representation $V$; and the eigenvalues of
the Casimir operator have the same property. As a matter of fact, 
our normalization $(6.3.6)$ and $(6.3.13)$ as well as the use of the operator 
${\cal A}^\delta$ have been done with $(6.3.27)$--$(6.3.29)$ in mind.
More precisely we have now
$$
\eqalign{
\nu_V=i\kappa_j\quad&\hbox{or}\quad k-\txt{1\over2},\cr
\varrho_V(n)={\Gamma(\txt{1\over2}+i\kappa_j)\over
2\pi^{{1\over2}+i\kappa_j}}\rho_j(n)\quad&\hbox{or}\quad
(-1)^k{\Gamma(2k)^{1\over2}\over 2^{2k}\pi^{k+{1\over2}}}
\rho_{j,k}(n),
}\eqno(6.3.30)
$$
according as $V$ belongs to either the unitary principal or
holomorphic/anti-holomorphic discrete series. We have, in place of
$(6.3.22)$,
$$
V=\bigoplus_{p=-\infty}^\infty {\Bbb C}\lambda_V^{(p)},\eqno(6.3.31)
$$
for any $V$, where an obvious convention is in force when $V$ is
not in the unitary principal series.
\medskip
As to the action of Hecke operators, the definitions
$(3.1.3)$ and $(3.1.19)$ can be translated into
$$
\eqalignno{
T(n)f(\r{g})&=n^{-{1\over2}}
\sum_{\tau\in \varGamma\backslash{\Bbb M}(n)}
f(n^{-{1\over2}}\tau\r{g})\cr
&=n^{-{1\over2}}\sum_{d|n}\sum_{b=1}^d f(\r{n}[b/d]\r{a}[n/d^2]\r{g}),
&(6.3.32)
}
$$
with ${\Bbb M}(n)$ as in Section 3.1. We have, in place of
$(3.1.13)$ and $(3.1.20)$,
$$
T(n)\lambda_V^{(p)}=t_V(n)\lambda_V^{(p)},\eqno(6.3.33)
$$
with $t_V(n)=t_j(n)$ or $t_{j,k}(n)$. This is due to $(6.3.5)$ and $(6.3.12)$
as well as to that $T(n)$ is defined with left translations in $(6.3.33)$;
that is, $T(n)$ commutes with the action of $\f g$, and
each $V$ is Hecke invariant. We have, for any integer 
$n$,
$$
\varrho_V(n)=\cases{\epsilon_V^{{1\over2}
(1-\sgn(n))}\varrho_V(1)t_V(|n|)& unitary principal series,\cr
{1\over2}(1+\sgn(n))\varrho_V(1)t_V(n)& discrete series,
}\eqno(6.3.34)
$$
where $\epsilon_V=\epsilon_j$ is as in $(3.1.15)$. The bound $(4.4.4)$
implies that uniformly for all $V$
$$
t_V(n)\ll n^{{1\over4}+\delta},\eqno(6.3.35)
$$
with the implied constant depending only on an arbitrary constant $\delta>0$.
\par
Finally, the automorphic $L$-functions $L_j$, $L_{j,k}$,
and their corresponding Hecke series $H_j$, $H_{j,k}$ which are 
introduced in Section 3.2 are
replaced by
$$
L_V(s)=\sum_{n=1}^\infty \varrho_V(n)n^{-s},\quad H_V(s)=\sum_{n=1}^\infty
t_V(n)n^{-s},\eqno(6.3.36)
$$
together with the normalization $(6.3.30)$. Also the definition
$(3.2.3)$ of the Rankin $L$-function is now extended by
$$
L_{V\otimes V'}(s)=\zeta(2s)\sum_{n=1}^\infty\varrho_V(n)
\overline{\varrho_{V'}(n)}n^{-s},
\eqno(6.3.37)
$$
for any pair $V,V'$ of irreducible representations.
\vskip 1cm
\centerline{\bf 6.4 Realization of representations}
\medskip
\medskip
\noindent
We shall try to investigate the structure of individual
subspaces $V$ listed in $(6.3.22)$ by means of a realization
of representations. This will lead us, in particular, to a
geometrical understanding of those exotic Bessel transforms $(2.3.17)$,
$(2.4.8)$, $(2.5.7)$, and $(2.5.15)$ which are involved in the
sum formulas of Kloosterman sums. More precisely, 
we shall find that those integral transforms are closely 
related to the action of the Weyl element 
$\r{w}=\left[{\atop-1}{1\atop}\right]$ in each 
subspace $V$. One should note here that the Bruhat decomposition
$$
\r{G}=\r{NA}\cup\r{NwNA}\eqno(6.4.1)
$$
holds, as we have 
$$
\left[\matrix{a&b\cr c& d}\right]=
\cases{\hfil\r{n}[ab]\r{a}[a^2]& if $c=0$,\cr
\r{n}[a/c]\r{w}\r{n}[cd]\r{a}[c^2]& if $c\ne0$.
}\eqno(6.4.2)
$$
Thus the realization of the actions of $\r{w}$, $\r{n}[x]$, $\r{a}[y]$ 
in each $V$ is a fundamental issue, and an answer to it is given in
Lemma 6.1 below.
\medskip
In the formulas $(6.3.9)$, $(6.3.16)$, and $(6.3.19)$ we see a
correspondence between the base elements of $V$ and the simple function
$\phi_p(\r{g},\nu_V)$ in terms of the operator ${\cal A}^\pm$ and
Fourier expansions of cusp forms over $\varGamma\backslash{\eusm H}$. 
The metrical properties $(6.3.10)$ and $(6.3.17)$ are 
in fact solely due to $(6.2.25)$, $(6.3.2)$ and 
the original orthonormality of the base system of those cusp forms.
As each subspace $V$ is generated 
from a particular cusp form and the
actions of the Maass operators upon it, 
it is naturally surmised that the structure of $V$ could
be described only with the behaviour of ${\cal A}^\pm$ in the
space spanned by $\phi_p$, without recoursing
to the nature of the cusp form generating $V$. In this way, we reach
the notion of {\it models\/}. What we are about to develop is, 
in fact, such a model.
\medskip
Thus, returning to $(6.3.1)$, we consider first the validity of the equation
$$
{\cal A}^\delta\!\phi(\r{g},\nu)
=\sum_{p=-\infty}^\infty c_p{\cal A}^\delta\!\phi_p(\r{g},\nu),
\eqno(6.4.3)
$$ 
with 
$$
\phi(\r{g},\nu)=\sum_{p=-\infty}^\infty c_p\phi_p(\r{g},\nu).\eqno(6.4.4)
$$
Here $\phi$ is assumed to be smooth, i.e., $c_p$ decays faster 
than any negative power of $|p|$ as
$p$ tends to infinity. Integration by parts gives
$$
\eqalignno{
&{\cal A}^\delta\phi_p(\r{a}[y],\nu)\cr
&= {y^{-{1\over2}-\nu}\over2\pi i}
\int_{-\infty}^\infty{((1+2\nu)\xi+2\delta pi)e(y\xi)
\over(\xi^2+1)^{{3\over2}+\nu}}
\left({\xi+i\over\xi-i}\right)^{\delta p} d\xi,&(6.4.5)
}
$$ 
which implies that $(6.4.3)$ holds for $\Re\nu>-{1\over2}$ via analytic
continuation with respect to $\nu$. Repeating the same procedure, we find
that $(6.4.3)$ holds for any $\nu\in\Bbb C$.
\medskip
Then, we define the Kirillov operator ${\eusm K}$ by
$$
{\eusm K}\phi(u)={\cal A}^{\sgn(u)}\phi(\r{a}[|u|]),
\quad u\in{\Bbb R}^\times,
\eqno(6.4.6)
$$ 
where ${\Bbb R}^\times={\Bbb R}\!\setminus\!\{0\}$; 
hereafter we shall often omit the parameter $\nu$ to avoid otherwise heavy
notations. This concept will play a crucial r\^ole 
in the sequel, for
it gives a realization of each $V$ in terms of a fairy ordinary
function space, that is, $L^2({\Bbb R}^\times,\pi^{-1}d^\times)$, 
$d^\times\! u=du/|u|$, as is made explicit in the following three lemmas.
\medskip
\medskip
\noindent 
{\bf Lemma 6.1} {\it Let $\phi$ be smooth as in
$(6.4.4)$. We have, with the right translation $\omega$,
$$
{\eusm K}\omega({\r{n}[x]})\phi(u)=e(ux){\eusm K}\phi(u),\quad
{\eusm K}\omega({\r{a}[y]})\phi(u)={\eusm K}\phi(uy).\eqno(6.4.7)
$$ 
Also, if $|\Re\nu|<{1\over2}$, then
$$
{\eusm K}\omega(\r{w})\phi(u)=\int_{{\Bbb R}^\times}j_\nu(u\lambda){\eusm K}
\phi(\lambda)d^\times\!\lambda.\eqno(6.4.8)
$$ 
Here
$$ 
j_\nu(u)=\pi{\sqrt{|u|}\over\sin\pi\nu}\left(J_{-2\nu}^{\sgn(u)}
(4\pi\sqrt{|u|})-J^{\sgn(u)}_{2\nu}(4\pi\sqrt{|u|})\right),
\eqno(6.4.9)
$$ 
where $J^+_\nu=J_\nu$ and $J^-_\nu=I_\nu$ with the ordinary
notation for Bessel functions. 
}
\medskip
\noindent 
{\it Proof\/}. Since $(6.4.7)$ is immediate, we deal with $(6.4.8)$ only.
For this sake we consider the Mellin transform
$$
\Gamma_p(s)=\Gamma_p(s,\nu)=\int_0^\infty
{\cal A}^+\!\phi_p(\r{a}[y])y^{s-{3\over2}}dy.\eqno(6.4.10)
$$ 
We shall show that $\Gamma_p(s)$ continues meromorphically to ${\Bbb C}$, and
satisfies the {\it local\/} functional equation
$$
\eqalignno{ (-1)^p\Gamma_p(s) =&2^{1-2s}\pi^{-2s}\Gamma(s+\nu)
\Gamma(s-\nu)\cr &\times\left(\cos\pi
s\,\Gamma_p(1-s)+\cos\pi\nu\,
\Gamma_{-p}(1-s)\right), &(6.4.11)
}
$$ 
provided $\Re\nu>-{1\over2}$.
In fact, by the first line of $(6.3.4)$,
we have, for $\Re s>\Re\nu>0$,
$$
\Gamma_p(s)=\int_0^\infty
y^{s-\nu-1}
\int_{\Im\xi={1\over2}}{e(y\xi)\over(\xi^2+1)^{\nu+{1\over2}}}
\left({\xi+i\over\xi-i}\right)^pd\xi dy.
$$  
Assuming temporarily that $0<\Re\nu<{1\over4}<\Re s<{1\over2}$,
we exchange the order of integration, and compute the resulting
inner integral. We find that 
$$
\eqalignno{
\Gamma_p(s)
&=(2\pi)^{\nu-s}
\Gamma(s-\nu)\cr
&\times\Big[\exp(-\txt{1\over2}\pi i(s-\nu))\r{L}_p(s)
+\exp(\txt{1\over2}\pi i(s-\nu))\r{L}_{-p}(s)\Big],&(6.4.12)
}  
$$  
with
$$ 
\r{L}_p(s)=\int_0^\infty{\xi^{-s+\nu}
\over(\xi^2+1)^{\nu+{1\over2}}}
\left({\xi+i\over\xi-i}\right)^pd\xi.\eqno(6.4.13)
$$
By analytic continuation, the expression $(6.4.12)$ holds if $-\Re\nu<\Re
s<1+\Re\nu$.  Under this condition,  we observe that the change of 
variable $\xi\to\xi^{-1}$ gives
$\r{L}_p(s)=(-1)^p\r{L}_{-p}(1-s)$. Then a rearrangement gives $(6.4.11)$ and
the meromorphic continuation of $\Gamma_p(s)$ follows via analytic continuation.
\par
We are now going to show that $(6.4.8)$ with $\phi=\phi_p$
and $(6.4.11)$ are in fact a Mellin pair; that is, the Mellin
inversion of $(6.4.11)$ yields a special case of $(6.4.8)$.
To this end, we note first that if $|\Re\nu|<\Re s<{1\over4}$, then
$$
\int_0^\infty j_\nu(u)u^{s-{3\over2}}du=2^{1-2s}\pi^{-2s}
\cos(\pi s)\Gamma(s+\nu)\Gamma(s-\nu),\eqno(6.4.14)
$$
and that if $|\Re\nu|<\Re s$, then
$$
\int_{-\infty}^0 j_\nu(u)|u|^{s-{3\over2}}du=2^{1-2s}\pi^{-2s}
\cos(\pi\nu)\Gamma(s+\nu)\Gamma(s-\nu).\eqno(6.4.15)
$$
The former is a consequence of $(4.4.11)$ and the latter of
$(3.2.38)$ with $(1.1.27)$; both integrals are absolutely
convergent. One could then appeal to the Parseval formula for $L^2$-pairs
of Mellin transforms. Here we develop instead a direct reasoning.
Thus, by the last two formulas we transform $(6.4.11)$ into
$$
(-1)^p\Gamma_p(s)=\int_{{\Bbb R}^\times}j_\nu(u)|u|^{s-{1\over2}}
\Gamma_{\sgn(u)p}(1-s)d^\times\!u.
$$
We replace $p$ by $\sgn(v)p$,
multiply both sides by the factor $|v|^{{1\over2}-s}/2\pi i$, and integrate
over the vertical line $\Re s=\beta$, with $|\Re\nu|<\beta<{1\over4}$. 
We have
$$
\eqalignno{
&{(-1)^p\over2\pi i}\int_{(\beta)}\Gamma_{\sgn(v)p}(s)|v|^{{1\over2}-s}ds\cr
&=\int_{{\Bbb R}^\times}j_\nu(u)
\Big\{{1\over2\pi i}\int_{(\beta)}\Gamma_{\sgn(uv)p}(1-s)
|u/v|^{s-{1\over2}}ds\Big\}
d^\times\!u. &(6.4.16)
}
$$
The absolute convergence that is needed to verify the exchange of the
order integration is due to the exponential decay of $\Gamma_p(s)$ which
can be confirmed by turning the contour in $(6.4.13)$
through a small angle appropriately. The left side of $(6.4.16)$ is equal
to 
$$
(-1)^p{\cal A}^+\!\phi_{\sgn(v)p}(\r{a}[|v|])
={\cal A}^{\sgn(v)}\!\phi_p(\r{a}[|v|]\r{w})={\eusm K}\omega(\r{w})\phi_p(v)
$$ 
in view of the first line of $(6.3.4)$; also, the inner-integral to 
$$
{\cal A}^+\!\phi_{\sgn(uv)p}(\r{a}[|u/v|])={\eusm K}\phi_p
(u/v),
$$
which yields $(6.4.8)$ in the case of $\phi=\phi_p$ with $|\Re\nu|<{1\over4}$.
To widen this range of $\nu$, we remark that
we have, for $\Re\nu>-{1\over2}$, $0<y\le1$,
$$
{\cal A}^\delta\phi_p(\r{a}[y])\ll
(|p|+|\nu|+1)y^{{1\over2}-|\Re\nu|}|\log y|,\eqno(6.4.17)
$$ 
and, for $\Re\nu>-{1\over2}$, $y\ge1$,
$$
{\cal A}^\delta\phi_p(\r{a}[y])\ll
(|p|+|\nu|+1)y^{-{1\over2}-\Re\nu}\exp\left(-{y\over  
|\nu|+|p|+1}\right).
\eqno(6.4.18)
$$ 
The implied constants in both bounds are absolute. In fact,
the first line of $(6.3.4)$ gives
$$
\eqalignno{
{\cal A}^\delta\phi_p(\r{a}[y])
&={\cal A}^\delta\phi_0(\r{a}[y])+y^{{1\over2}-\nu}\int_{-\infty}^\infty
{e(y\xi)\over (\xi^2+1)^{{1\over2}+\nu}}
\left(\left({\xi+i\over\xi-i}\right)^{\delta p}-1\right)
d\xi\cr
&={2\pi^{{1\over2}+\nu}\over\Gamma({1\over2}+\nu)}y^{1\over2}
K_\nu(2\pi y) +O\left(y^{{1\over2}-\Re\nu}(|p|+1)\right), &(6.4.19)
}
$$ 
provided $\Re\nu>-{1\over2}$, and $(6.4.17)$ follows. As to $(6.4.18)$,
it suffices to shift the contour in $(6.4.5)$ to 
$\Im\xi=(|\nu|+|p|+1)^{-1}$. Via these bounds we get
the desired analytic continuation to $|\Re\nu|<{1\over2}$. The assertion
$(6.4.8)$ with a general smooth $\phi$ is now immediate. 
This ends the proof of the lemma.
\medskip
Next, we shall show that the Kirillov operator is in fact a unitary map of
a simple nature:
\medskip
\medskip
\noindent {\bf Lemma 6.2.}\quad{\it We assume that $\nu\in i{\Bbb R}$, 
and introduce the Hilbert space
$$ 
U_\nu=\bigoplus_{p=-\infty}^\infty
{\Bbb C}\phi_p,\quad
\phi_p(\r{g})=\phi_p(\r{g};\nu),
\eqno(6.4.20)
$$ 
equipped with the ordinary norm
$$
\Vert\phi\Vert_{U_\nu}^2=\sum_{p=-\infty}^\infty
|c_p|^2,\quad\phi= \sum_{p=-\infty}^\infty c_p\phi_p.\eqno(6.4.21)
$$ 
Then the operator ${\eusm K}$ is a unitary map from $U_\nu$ onto
$L^2({\Bbb R}^\times,\pi^{-1}d^\times)$. In particular,
$\omega$ and ${\eusm K}\omega{\eusm K}^{-1}$ are equivalent 
unitary representations of $\r{G}$ in $U_\nu$ and
$L^2({\Bbb R}^\times,\pi^{-1}d^\times)$, respectively.
\/}
\medskip
\noindent 
{\it Proof\/}.  We shall first treat the second assertion, while
assuming the validity of the first. The
unitarity of ${\eusm K}\omega(\r{n}[x]){\eusm K}^{-1}$ and 
${\eusm K}\omega(\r{a}[y]){\eusm K}^{-1}$ 
with respect to $L^2({\Bbb R}^\times,\pi^{-1}d^\times)$ is obvious
from $(6.4.7)$. Also, the unitarity of $\omega(\r{k}[\theta])$ 
on $U_\nu$ is fairy obvious. Hence the assertion follows.
\par
Let us prove the unitarity of $\eusm K$. We shall 
employ an explicit reasoning. Thus, by $(6.3.4)$ and
$(6.4.6)$, we have, for
$\nu\in i{\Bbb R}$ and $p,q\in{\Bbb Z}$,
$$
\eqalignno{
&{1\over\pi}\int_{{\Bbb R}^\times}{\eusm K}\phi_p(u)
\overline{{\eusm K}\phi_q(u)}d^\times\!u\cr
&={(-1)^{p+q}\over\Gamma(p+\nu+{1\over2})
\Gamma(q-\nu+{1\over2})}\int_0^\infty W_{p,\nu}(y)
W_{q,\nu}(y){dy\over y}\cr
&+{(-1)^{p+q}\over\Gamma(-p+\nu+{1\over2})
\Gamma(-q-\nu+{1\over2})}\int_0^\infty W_{-p,\nu}(y)
W_{-q,\nu}(y){dy\over y},&(6.4.22)
}
$$
where we have used the fact that $W_{p,\nu}(y)$ is real because
of $(3.2.34)$. To evaluate these integrals, we appeal to the
following formula: For any
$\alpha,\beta\in{\Bbb C}$ and
$|\Re\mu|<{1\over2}$, it holds that
$$
\eqalignno{
&\int_0^\infty W_{\alpha,\mu}(y) W_{\beta,\mu}(y){dy\over y}=
{\pi\over(\alpha-\beta)\sin(2\pi\mu)}\cr
&\times\left[{1\over\Gamma({1\over2}-\alpha+\mu)\Gamma({1\over2}
-\beta-\mu)}-{1\over\Gamma({1\over2}-\alpha-\mu)\Gamma({1\over2}
-\beta+\mu)}\right],\quad&(6.4.23)
}
$$ 
together with
$$
\eqalignno{
&\int_0^\infty (W_{\alpha,\mu}(y))^2{dy\over y}=
{\pi\over\sin(2\pi\mu)}\cr
&\times{1\over\Gamma({1\over2}-\alpha+\mu)\Gamma({1\over2}
-\alpha-\mu)}\Big[{\Gamma'\over\Gamma}(\txt{1\over2}-\alpha+\mu)
-{\Gamma'\over\Gamma}(\txt{1\over2}-\alpha-\mu)\Big].\qquad&(6.4.24)
}
$$ 
By $(6.4.23)$ we see readily that the right side of $(6.4.22)$ vanishes whenever
$p\ne q$; and by $(6.4.24)$ it is equal to $1$ if $p=q$. Hence
$$
{1\over\pi}\int_{{\Bbb R}^\times}{\eusm K}\phi_p(u)
\overline{{\eusm K}\phi_q(u)}d^\times\!u=\delta_{p,q},\eqno(6.4.25)
$$
which is equivalent to the unitarity of $\eusm K$.
\par
To show $(6.4.23)$, we use the Whittaker differential equation 
$(3.2.32)$. We have
$$
\eqalignno{ 
&-\alpha\int_0^\infty
W_{\alpha,\mu}(y)W_{\beta,\mu}(y){dy\over y}\cr
&=\lim_{\varepsilon\to0^+}
\int_\varepsilon^\infty\Big[\left({d\over dy}
\right)^2-{1\over4}+{{1\over4}-\mu^2\over y^2}\Big]
W_{\alpha,\mu}(y)W_{\beta,\mu}(y){dy\over y}\cr
&=\lim_{\varepsilon\to0^+}\left[-W_{\alpha,\nu}'(\varepsilon)
W_{\beta,\mu}(\varepsilon)+W_{\alpha,\mu}
(\varepsilon)W_{\beta,\mu}' (\varepsilon)\right]\cr
&\hskip1.5cm-\beta\int_0^\infty
W_{\alpha,\mu}(y)W_{\beta,\mu}(y){dy\over y}.&(6.4.26)
}
$$ 
To compute the last limit we invoke that near the origin 
$$
\eqalignno{
W_{\alpha,\mu}(y)&=\!\Big({\Gamma(-2\mu)\over\Gamma({1\over2}-\alpha-\mu)}
y^{\mu+{1\over2}}+{\Gamma(2\mu)\over\Gamma({1\over2}-\alpha+\mu)}
y^{-\mu+{1\over2}}\Big)(1+O(y)),\cr
W'_{\alpha,\mu}(y)&=\!\Big({(\mu+{1\over2})
\Gamma(-2\mu)\over\Gamma({1\over2}-\alpha-\mu)}
y^{\mu-{1\over2}}-{(\mu-{1\over2})\Gamma(2\mu)
\over\Gamma({1\over2}-\alpha+\mu)}
y^{-\mu-{1\over2}}\Big)(1+O(y)),\qquad\qquad&(6.4.27)
}
$$
where the implied constant is bounded as far as $|\Re\mu|<{1\over2}$.
After a rearrangement we get $(6.4.23)$. 
\par
It remains for us to show the surjectivity of $\eusm K$. 
Thus, let $\nu\in i{\Bbb R}$ and assume that a smooth   function
$g$, compactly supported on ${\Bbb R}^\times$, is orthogonal to all
${\eusm K}\phi_p$. Multiply $(6.4.5)$ by $g$ and integrate, change the
order of integration, and undo the integration by parts with respect
to the outer integral. We have
$$
\leqalignno{ 0&=\int_{{\Bbb R}^\times}g(u)\overline{{\eusm K}\phi_p(u)}
d^\times\!u\cr &=\int_{-\infty}^\infty
{1\over(\xi^2+1)^{{1\over2}+\nu}}
\left({\xi+i\over\xi-i}\right)^p\int_{-\infty}^\infty 
g(u) |u|^{-{1\over2}+\nu}e(-u\xi)du\,d\xi. }
$$ 
Then we note that the system
$\left\{((\xi+i)/(\xi-i))^p:\,p\in{\Bbb Z}\right\}$ is complete
orthonormal in the space  
$L^2\!\left({\Bbb R},(\pi(\xi^2+1))^{-1}d\xi\right)$. Hence the Fourier
transform of $g(u)|u|^{-{1\over2}+\nu}$ vanishes identically, whence
the assertion. This ends the proof of the lemma.
\medskip 
In passing, we make a remark on the complementary series, i.e., the
situation with $-{1\over2}<\nu<{1\over2}$, although such a 
representation of $\r{G}$ does not occur in
$L^2(\varGamma\backslash\r{G})$.  It is easy to see that
Lemma 6.1 remains valid. The definition $(6.4.20)$ is the same, 
but $(6.4.21)$ has to be replaced by the  norm
$$
\pi^{\nu}\Big(\sum_{p=-\infty}^\infty{\Gamma(p+{1\over2}-\nu)
\over\Gamma(p+{1\over2}+\nu)} |c_p|^2\Big)^{1\over2}\,.\eqno(6.4.28)
$$ 
With this, the above proof extends readily, and Lemma 6.2 holds for
these $\nu$ as well.
\medskip 
On the other hand, in dealing with the holomorphic discrete
series,  
$(6.4.20)$ has to be replaced by the Hilbert space
$$ 
D_k=\bigoplus_{p=k}^\infty{\Bbb C}\phi_p,\quad \phi_p(\r{g})=
\phi_p\left(\r{g};k-\txt{1\over2}\right),\eqno(6.4.29)
$$ 
with an integer $k\ge1$, which is equipped with the norm
$$
\Vert\phi\Vert_{D_k}=\pi^{k-{1\over2}}\Big(\sum_{p=k}^\infty
{\Gamma(p-k+1)\over\Gamma(p+k)}|c_p|^2\Big)^{1\over2},
\quad \phi=\sum_{p=k}^\infty c_p\phi_p.\eqno(6.4.30)
$$ 
Since ${\cal A}^-$ annihilates $D_k$, we are concerned with
${\cal A}^+$ only. The expression $(6.3.4)$, $\delta=+$, holds without
changes. With this, the operator ${\eusm K}$ is defined by $(6.4.6)$ again;
it should be noted that
$$
j_{k-{1\over2}}(u)=\cases{\hfil 0, & if $u<0$,\cr
2\pi(-1)^k\sqrt{u} J_{2k-1}(4\pi\sqrt{u}),& if $
u>0$.
}\eqno(6.4.31)
$$
\medskip
\noindent 
{\bf Lemma 6.3.} {\it The operator ${\eusm K}$ is a unitary map
of $D_k$ onto $L^2((0,\infty),\pi^{-1}d^\times)$. Also, for any
smooth $\phi\in D_k$, we have $(6.4.7)$ and $(6.4.8)$ with $\nu=k-{1\over2}$.
With these changes the second assertion of the previous lemma holds.
The analogue for the anti-holomorphic discrete series is obtained by  
applying the involution $\r{g=nak}\mapsto
\r{n}^{-1}\r{a}\r{k}^{-1}$. 
}
\medskip
\noindent 
{\it Proof.\/} The third assertion is immediate. As to the
unitarity of ${\eusm K}$, it is proved with a minor change of the
above argument. In fact, the Whittaker function
$W_{p,k-{1\over2}}(u)$ ($p\ge k$) is a product of
$u^k\exp(-{1\over2}u)$ and a polynomial on $u$ of degree $p-k$, as
$(6.3.4)$ implies. Thus the proof of $(6.4.23)$ can be carried
out also for the product
$W_{p,k-{1\over2}}(u)W_{q,k-{1\over2}}(u)$ with integers
$p,\,q$, although the condition on $\Re\mu$ there is violated. The result
is equal to the limit of $(6.4.23)$ as $(\alpha,\beta,\nu)$ tends to
$\left(p,q,k-{1\over2}\right)$.  About the surjectivity, we
argue as follows: Let $g$ be smooth and compactly supported on
$(0,\infty)$. If $g$ is orthogonal to all
${\eusm K}\phi_p$, $p\ge k$, then we have, by the remark just made on
$W_{p,k-{1\over2}}(u)$,
$$
\int_0^\infty g(u)\exp(-2\pi u)u^{p-1}du=0,\quad p\ge k.
\eqno(6.4.32)
$$ 
This implies that the Fourier transform of
$g(u)\exp(-2\pi u)u^{k-1}$ vanishes identically; in fact it
suffices to expand the additive character into a power series and
integrate termwise. Hence $g\equiv0$. On noting $(6.4.31)$,
the counterpart of $(6.4.8)$, with $\phi=\phi_p$, 
$\nu=k-{1\over2}$, can be proved in much the same way as before.
The extension to any smooth $\phi$ is immediate via $(4.4.24)$ and
$$
{\eusm K}\phi_p(u)={\cal A}^+\!\phi_p(\r{a}[u])
\ll\min(u,|p|+1)u^{-k},\quad u>0,\eqno(6.4.33)
$$ 
which comes from $(6.3.4)$ and $(6.4.5)$. This ends the
proof of the lemma.
\medskip
Here we summarize our discussion in the present section: 
Let $V$ be a subspace listed in $(6.3.22)$; we assume 
for instance that $\nu_V\in i{\Bbb R}$.
We put ${\eusm L}(\lambda^{(p)}_V)=\phi_p(\cdot,\nu_V)$,  
and extend $\eusm L$ to the whole of $V$
in an obvious manner. Then ${\eusm L}$ is an isometry mapping $V$
onto $U_{\nu_V}$, where the metric of the latter is defined in
Lemma 6.2. We put
$$
\omega_V=({\eusm KL})\omega({\eusm KL})^{-1},\eqno(6.4.34)
$$
with the right translation $\omega$ acting in $V$. 
Then $\omega_V$ is a representation of $\r{G}$ 
in the Hilbert space $L^2({\Bbb R}^\times,\pi^{-1}d^\times)$.
In this way, the representation $V$
is realized in terms of $L^2({\Bbb R}^\times\!,\pi^{-1}d^\times)$. 
On the other hand, $\omega_\nu={\eusm K}\omega{\eusm K}^{-1}$
with $\omega$ acting in $U_\nu$
is also a representation of $\r{G}$ in $L^2({\Bbb
R}^\times,\pi^{-1}d^\times)$, for any $\nu\in i{\Bbb R}$. 
The explicit description of the actions of $\r{G}$ under
$\omega_\nu$ is given in Lemma 6.1. The extension to the discrete
series of representations is given in Lemma 6.3.
\medskip
Now, as we have promised, we shall give a proof of Theorem 6.2. We may assume
that $V$ be such that $\nu_V\in i{\Bbb R}$, for other cases are in fact
easier. Naturally, it suffices to prove that $\omega_\nu$ is an
irreducible representation. Let $Y_1$ be an invariant 
subspace of $L^2({\Bbb R}^\times\!,\pi^{-1}d^\times)$, and $Y_2$ its
orthogonal complement, which is also invariant.
Let $f_j\in Y_j$ be arbitrary; hereafter, equalities are in the $L^2$-sense. 
Since $\omega_\nu(\r{n}[x])f_1\in Y_1$ for any real $x$, the orthogonality
of $Y_1$ and $Y_2$ implies that the Fourier
transform of $f_1(u)\overline{f_2(u)}/|u|$ vanishes identically, because of 
the first identity in $(6.4.7)$. That is, $f_1f_2=0$. 
Then, by the second identity in $(6.4.7)$ we see that $f_1(uy)f_2(u)=0$ for
any $y>0$. This means that $f_1(u)f_2(v)=0$ for $uv>0$. Consequently we may assume
without loss of generality that any $f\in Y_1$ is such that $f(u)=0$ 
for $u<0$. By $(6.4.8)$ we have, for $u<0$,
$$
0=\omega_\nu(\r{w})f(u)=\int_0^\infty
j_{\nu}(u\lambda)f(\lambda){d\lambda\over\lambda}.\eqno(6.4.35)
$$
Let $\tilde{f}(s)$ be the $L^2$-Mellin transform of $f$, which should
exist for $s\in i{\Bbb R}$. We observe that  by $(6.4.15)$
$$
\eqalign{
&\int_0^\infty j_{\nu}(u\lambda)\lambda^{s-2}d\lambda\cr
&=4\pi\cos(\pi\nu)(2\pi)^{-2s}\Gamma(s-\txt{1\over2}+\nu)
\Gamma(s-\txt{1\over2}-\nu)|u|^{1-s},
}
$$
for $\Re s>{1\over2}$ and $u<0$. Hence, by the Mellin-Parseval identity,
$(6.4.35)$ is equivalent to
$$
\int_{(1)}(2\pi)^{-2s}\Gamma(s-\txt{1\over2}+\nu)
\Gamma(s-\txt{1\over2}-\nu)\tilde{f}(1-s)|u|^{-s}ds=0,
$$
for any $u<0$. The integrand has to vanish, and $\tilde{f}(s)=0$ for
$s\in i{\Bbb R}$. Namely, we have $\Vert{f}\Vert=0$ with the norm
in $L^2({\Bbb R}^\times\!,\pi^{-1}d^\times)$. We end the proof of Theorem 6.2.
\vskip 1cm
\centerline{\bf 6.5 Revisits}
\medskip
\noindent
The aim of this section is to review the sum formulas of Kloosterman sums
and the explicit formula for the fourth moment of the
zeta-function in the light of automorphic representations. We shall, however,
restrict ourselves to the structural aspect of the new argument leading to
those formulas, as a fuller account would not be a help, rather a hindrance
to see the essentials.
Thus, for instance, convergence issues
will be ignored. We shall first discuss the sum formulas and then turn to
the zeta-function.
\medskip
Thus, let us reformulate the sum formulas:
For any non-zero integers $m,n$ and for appropriate weight
functions $f, \varphi$, we have
$$
\eqalignno{
&\sum_V{\raise4pt\hbox{${}\!^{\r{u}}$}}
\,\overline{\varrho_V(m)}\varrho_V(n)f(\nu_V)
+ {1\over4\pi i}\int_{(0)}{\sigma_{2r}(m)\sigma_{2r}(n)\over
(mn)^r\zeta(1+2r)\zeta(1-2r)}f(r)dr\cr
=&\delta_{mn}{i\over4\pi^2}\int_{(0)}r\tan(\pi r)f(r)dr
+\sum_{\ell=1}^\infty{1\over\ell}S(m,n;\ell)
{\bf A}\!^\delta\!f\Big({4\pi\over\ell}
\sqrt{|mn|}\Big),\quad&(6.5.1)
}
$$
as well as
$$
\eqalignno{
\sum_{\ell=1}^\infty&{1\over\ell}S(m,n;\ell)\varphi
\Big({4\pi\over\ell}\sqrt{|mn|}\Big)
=\sum_V\overline{\varrho_V(m)}\varrho_V(n){\bf B}^\delta\!\varphi(\nu_V)\cr
&+{1\over4\pi i}\int_{(0)}{\sigma_{2r}(m)\sigma_{2r}(n)\over
(mn)^r\zeta(1+2r)\zeta(1-2r)}{\bf B}^\delta\!\varphi(r)dr.&(6.5.2)
}
$$
Here $\delta=\sgn(mn)$ and $\sum^\r{u}$ indicates that the sum is
restricted to all irreducible representations in
the unitary principal series; in $(6.5.2)$
$V$ runs over all irreducible representations listed in $(6.3.22)$. Also,
$$
\eqalignno{
{\bf A}^\delta\!f(x)&={i\over4\pi}\int_{(0)}
{J^\delta_{-2\nu}(x)-J^{\delta}_\nu(x)\over\sin(\pi\nu)}\nu\tan(\pi\nu)
f(\nu)d\nu,
&(6.5.3)\cr
{\bf B}^\delta\!\varphi(\nu)&={2\pi}\int_0^\infty
{J^\delta_{-2\nu}(x)-J^{\delta}_\nu(x)\over\sin(\pi\nu)}\varphi(x){dx\over
x},&(6.5.4) }
$$
where $J^\pm_{2\nu}$ are as in $(6.4.9)$. The formula $(6.5.1)$ is
equivalent to Theorems 2.2 and 2.4; $(6.5.2)$ to Theorems 2.3.and 2.5;
also, $(6.5.3)$ to $(2.3.17)$ and $(2.5.7)$; and $(6.5.4)$ to $(2.4.8)$
and $(2.5.15)$. Note that the present $f,\varphi$
are not the same as in those theorems,
which is due partly to the renormalization $(6.3.30)$.
\medskip
Although it does not matter for practical purposes, 
our proofs of Theorems 2.2--2.5 or $(6.5.1)$--$(6.5.2)$ that are
developed in Chapter 2 are admittedly highly
technical.  Also, the  emergence of holomorphic cusp forms 
in the statement of Theorem 2.3 or the same
in $(6.5.2)$ remains  mysterious and unexplained. 
With what we have developed in the present chapter, 
one may infer that the former phenomenon  should be related
to the assertion $(6.4.8)$ as a minor modification of
the kernel $j_\nu$, {\it the
Bessel function of representations\/} under our specification with 
$\r{G}=\r{PSL}(2,{\Bbb R})$, appears in the Bessel transforms
${\bf A}^\delta,{\bf B}^\delta$. Also
the latter phenomenon should be related to the spectral decomposition $(6.3.20)$
with $(6.3.21)$--$(6.3.24)$, as there cusp forms of all types play
their r\^oles, without any notable discrimination among them,
through irreducible representations generated by them.
That is, the sum formulas should better be
captured as a consequence of Theorem 6.1 augmented by Lemmas
6.1--6.3, not as that of Theorem 1.1. Namely, we need to
devise an argument to prove the sum formulas in the framework
of the space $L^2(\varGamma\backslash\r{G})$.
\medskip
To this end, we consider a Poincar\'e series on $\r{G}$. 
Let the seed function $q(\r{g})$ be defined 
on the big Bruhat cell, i.e., the set $\r{NwNA}$ in $(6.3.1)$ 
in such a way that $q(\r{n}[x_1]\r{w}
\r{n}[x_2]\r{a}[u])=\exp(2\pi imx_1)
\eta(x_2)g(u)$ with an integer $m\ne0$ and smooth functions
$\eta,\, g$ of fast decay; and $q(\r{g})=0$ on the
small cell, i.e., $\r{NA}$. In particular we have
$q(\r{n}[\xi]\r{g})=\exp(2\pi im\xi)q(\r{g})$ for any real $\xi$. We then put
$$
{\eusm Q}(\r{g})=\sum_{\gamma\in\varGamma_\infty
\backslash\varGamma}q(\gamma\r{g}),\eqno(6.5.5)
$$
with $\gamma$ as in $(1.1.4)$. In much the same way as the derivation
of $(1.1.6)$, we have
$$
{\eusm Q}(\r{g})=q(\r{g})
+\sum_{n=-\infty}^\infty
\sum_{\ell=1}^\infty S(m,n;\ell)
\int_{-\infty}^\infty
e(-n\xi)q(\r{n}[x_1]\r{a}[y_1]\r{k}[\theta_1])d\xi, \eqno(6.5.6)
$$
where
$$
\eqalign{
x_1=-&{\xi+x\over\ell^2((\xi+x)^2+y^2)},
\quad y_1={y\over\ell^2((\xi+x)^2+y^2)},\cr
&\exp(i\theta_1)=\exp(i\theta){\xi-iy\over(\xi^2+y^2)^{1\over2}}.
}
$$
We observe that
$$
\r{n}[x_1]\r{a}[y_1]\r{k}[\theta_1]=\r{a}[\ell^{-2}]\r{w}\r{n}[\xi]\r{g},
\eqno(6.5.7)
$$
and that the last integral at $\r{g}=1$ is equal to
$$
\int_{-\infty}^\infty e(-n\xi)q
(\r{w}\r{n}[\ell^2\xi]\r{a}[\ell^2])d\xi=\ell^{-2}
\hat{\eta}(2\pi n/\ell^2)g(\ell^2),\eqno(6.5.8)
$$
with the Fourier transform $\hat{\eta}$.
\par
On the other hand, the projection of ${\eusm Q}(\r{g})$ to an irreducible
subspace $V$ in the unitary principal series is equal to
$$
\sum_{p=-\infty}^\infty\langle{{\eusm Q},\lambda_V^{(p)}}\rangle 
\lambda_V^{(p)}(\r{g})=\sum_{p=-\infty}^\infty\lambda_V^{(p)}(\r{g})
\int_{\varGamma_\infty\backslash\r{NwNA}}q(\r{h})
\overline{\lambda_V^{(p)}(\r{h})}d\r{h}.
$$
The value at $\r{g}=1$ of the $n$-th Fourier 
coefficient of this expression is equal to
$$
\eqalignno{
&{\overline{\varrho_V(m)}\varrho_V(n)\over
\sqrt{|mn|}}\sum_{p=-\infty}^\infty{\cal A}^{\sgn(n)}\phi_p(\r{a}[|n|])\cr
&\times\int_0^\infty\int_{-\infty}^\infty \eta(\xi)g(u)
\overline{{\cal A}^{\sgn(m)}\phi_p(\r{a}[|m|]\r{w}\r{n}[\xi]\r{a}[u])}d\xi 
{du\over \pi u},&(6.5.9)
}
$$
where $\phi_p=\phi_p(\cdot,\nu_V)$, and 
we have used the fact that the Jacobian of the
change of variables
$$
\r{n}[x]\r{a}[y]\r{k}[\theta]\mapsto \r{n}[x_1]\r{w}\r{n}[\xi]\r{a}[u]
$$
is equal to $y^2/u$, that is, $d\r{g}=dx_1d\xi du/(\pi u)$. We have
$$
\eqalignno{
{\cal A}^{\sgn(m)}&\phi_p(\r{a}[|m|]\r{w}\r{n}[\xi]\r{a}[u])
={\eusm K}
\omega(\r{w}\r{n}[\xi]\r{a}[u])\phi_p(m)\cr
&=\int_{{\Bbb R}^\times}j_{\nu_V}(m\lambda)
\omega(\r{n}[\xi]\r{a}[u]){\eusm K}\phi_p(\lambda)d^\times\!\lambda\cr
&=\int_{{\Bbb R}^\times}j_{\nu_V}(m\lambda)e(\xi\lambda)
{\eusm K}\phi_p(u\lambda)d^\times\!\lambda,&(6.5.10)
}
$$
in which the second line is due to $(6.4.7)$, and the third to $(6.4.8)$.
Thus the double integral in $(6.5.9)$ is equal to
$$
\int_{{\Bbb R}^\times}\Big(\int_0^\infty
\hat{\eta}(2\pi\lambda/u)g(u)j_{\nu_V}(m\lambda/u){du\over u}\Big)
\overline{{\eusm K}\phi_p(\lambda)}d^\times\!\lambda,
$$
where we have used that $j_\nu$ is real valued. Then, the sum in $(6.5.9)$ is
equal to
$$
\pi\int_0^\infty
\hat{\eta}(2\pi n/u)g(u)j_{\nu_V}(mn/u){du\over u},\eqno(6.5.11)
$$
for ${\cal A}^{\sgn(n)}\phi_p(\r{a}[|n|])={\eusm K}\phi_p(n)$, and
$\{{\eusm K}\phi_p: p\in{\Bbb Z}\}$ is a complete orthonormal system of the space
$L^2({\Bbb R}^\times\!,\pi^{-1}d^\times)$, according to Lemma 6.2.
\par
Collecting these, we see that the contribution of 
the irreducible representation $V$ to the spectral
expansion of the sum
$$
\sum_{\ell=1}^\infty{1\over\ell^2}S(m,n;\ell)\hat{\eta}(2\pi n/\ell^2)
g(\ell^2)
$$
is equal to
$$
{\overline{\varrho_V(m)}\varrho_V(n)\over\sqrt{|mn|}}
\int_0^\infty\hat{\eta}(2\pi n/u)g(u)j_{\nu_V}(mn/u){du\over u}.
$$
Assuming that $m,n>0$,
we put $\hat{\eta}(x)=\tilde{\eta}(mx)$,
$g(u)=\tilde{g}(mn/u)$, and further 
$\tilde{\eta}(2\pi u)\tilde{g} (u)=\varphi(4\pi\sqrt{u})/(4\pi\sqrt{u})$. 
Then, we recover the integral transform ${\bf B}^+$ defined by $(6.5.4)$
and the relevant part of $(6.5.2)$. The contribution 
of the representations in the
discrete series and that of the continuous spectrum are treated
fairly analogously. In this way we have obtained a proof of Theorem 2.3,
i.e., $(6.5.2)$ with $mn>0$ via
the theory of automorphic representations, as far as we restrict
ourselves to the present choice of the weight function $\varphi$. Also,
the case with $mn<0$ can be treated similarly. The transforms
${\bf B}^\pm$ have turned out indeed to be equivalent to $(6.4.8)$ 
and its relevant statement given in Lemma 6.3.
\par
The mechanism can be summarized as this: The Fourier expansion of
Poincar\'e series with respect to the left action of $\r{N}$ takes us
to the notion of the operator ${\cal A}^\delta$ as $(6.5.7)$ dictates, 
and to the big Bruhat cell. The former demands
a theory of representations expressed in terms of ${\cal A}^\delta$, and
this leads us to the theory of the Kirillov model as 
we have given in Lemmas 6.1--6.3. The big Bruhat cell is of course
characterized by the presence of the Weyl element $\r{w}$, and its
action has to be {\it realized\/} if any practical application
of the harmonic analysis has to be performed. There naturally
emerges the Bessel function $j_\nu$  of representations.
On the other hand it now becomes expedient for us to work with
functions defined in the big Bruhat cell as $(6.5.8)$ shows clearly.
With this, the rest of the procedure is quite plain as 
$(6.5.9)$--$(6.5.11)$ is simply a logical rearrangement, although
it is true that the inversion argument at $(6.5.11)$
of a Fourier type is of some interest.
\medskip
The above discussion is, however, highly formal. There are a few missing points.
One is the treatment of the convergence issue; and the other is
the expansion of the space of weight functions so that the full statement
of Theorems 2.3 and 2.5 be recovered. These are, however, technical
issues, and could be regarded as being outside our present aim.   
A more essential problem than them is the derivation of Theorems 2.2 and
2.4 from 2.3 and 2.5, respectively; that is, the Spectral--Kloosterman
sum formula $(6.5.1)$ is to be derived from the Kloosterman--Spectral sum
formula $(6.5.2)$, the direction of which is exactly opposite to the reasoning 
in Sections 2.4 and 2.5. Its solution is naturally a logical necessity as far
as we proceed as in the present chapter. Our answer to this 
is that for both $\delta=\pm1$
$$
\hbox{\it The sum formulas $(6.5.1)$ and $(6.5.2)$ are 
equivalent to each other.} \eqno(6.5.12)
$$
We skip the proof, for it would be a digression too long.
\medskip
We turn to a review of Theorem 4.2. Thus, we should consider rather
$(4.2.5)$ than the theorem itself. The exploitation of the view point provided
there  has been the main motivation for us to develop an account on automorphic
representations. Since we are now working with matrices in projective
sense, $(4.2.5)$ should be reformulated as
$$
\eqalignno{
\sum_{n=1}^\infty n^{-z-{1\over2}}&\sum_{d|n}\sum_{b=1}^d
{\eusm P}_F(\r{n}[b/d]\r{a}[n/d^2]\r{g}),&(6.5.13)\cr
&{\eusm P}_F(\r{g})=\sum_{\gamma\in\varGamma}F(\gamma\r{g}),&(6.5.14)
}
$$
where the definition  $(6.3.32)$ of Hecke operators is taken into
account, and the convergent factor $n^{-z}$ is inserted, with $\Re z$
being sufficiently large. We are now to carry
out the computation of the  spectral decomposition of the Poincar\'e series
${\eusm P}_F(\r{g})$ via Theorem 6.1 and Lemmas 6.1--6.3; 
note that our interest is in fact
in the special value  at $\r{g}=1$.
Our argument is again formal. This time we
take into account the action of Hecke operators, i.e., the
assumption $(6.3.32)$--$(6.3.34)$.
\par
Thus, the value at $\r{g}=1$ of the projection of $(6.5.13)$ 
to an irreducible subspace $V$ in the unitary principal series is
equal to
$$
H_V(z)\sum_{p=-\infty}^\infty\langle 
{\eusm P}_F,\,\lambda_V^{(p)}\rangle\lambda_V^{(p)}(1),\eqno(6.5.15)
$$
where $H_V$ is as in $(6.3.35)$. This sum is 
$$
\eqalignno{
&\sum_{p=-\infty}^\infty\lambda_V^{(p)}(1)\int_G
F(\r{g})\overline{\lambda_V^{(p)}(\r{g})}d\r{g}\cr
=&\overline{\varrho_V(1)}\sum_{p=-\infty}^\infty\lambda_V^{(p)}(1)
\sum_{m=1}^\infty{t_V(m)\over\sqrt{m}}
\left(\Phi_p^++\epsilon_V\Phi_p^-\right)F(m,\nu_V)\cr
=&|\varrho_V(1)|^2\sum_{m=1}^\infty
\sum_{n=1}^\infty{t_V(m)t_V(n)\over\sqrt{mn}}\cr 
&\times\left(\e{B}^{(+,+)}+\e{B}^{(-,-)}+\epsilon_V\e{B}^{(+,-)}
+\epsilon_V\e{B}^{(-,+)}\right)F\left(\r{a}[n];m,\nu_V\right),
\qquad &(6.5.16)
}
$$ 
where
$$
\e{B}^{(\delta_1,\delta_2)}
F(\r{g};m,\nu_V)=\sum_{p=-\infty}^\infty\Phi_p^{\delta_1}F(m,\nu_V)\,
\e{A}^{\delta_2}\phi_p(\r{g}), \eqno(6.5.17)
$$ 
with
$$
\Phi_p^\delta F(m,\nu)=\int_G
F(\r{g})\overline{\e{A}^\delta\phi_p(\r{a}[m]\r{g})}d\r{g},
\quad \phi_p(\r{g})=\phi_p(\r{g},\nu).\eqno(6.5.18)
$$ 
We have, in terms of the Kirillov operator,
$$
\e{B}^{(\delta_1,\delta_2)}F(\r{a}[n];m,\nu)
=\sum_{p=-\infty}^\infty\Phi_p^{\delta_1}F(m,\nu)\,
{\eusm K}\phi_p(\delta_2 n).\eqno(6.5.19)
$$
We then proceed just as in $(6.5.9)$--$(6.5.11)$.
Since the integral in $(6.5.18)$ can be restricted to the big Bruhat cell, we
perform the change of variables accordingly. We have,
with $\r{h}=\r{n}[x_1]\r{w}\r{n}[x_2]\r{a}[u]$,
$$
\eqalignno{
\Phi_p^{\delta}F(m,\nu)&=\int_{\r{NwNA}}F(\r{h})
\overline{\omega(\r{h})\e{K}\phi_p(\delta m)}
dx_1dx_2{du\over\pi u}\cr
&=\int_{{\Bbb R}^\times}\Big(\int_{\r{NwNA}} 
F(\r{h})j_\nu(\delta m\lambda/u)\cr
&\qquad\times e(-\delta mx_1u-\lambda x_2/u)
dx_1dx_2{du\over\pi u}\Big)
\overline{\e{K}\phi_p(\lambda)}d^\times\!\lambda,
\quad\qquad&(6.5.20)
}
$$ 
where we have applied Lemma 6.1. Inserting this into $(6.5.19)$, we get,
via Lemma 6.2, that
$$
\eqalignno{ 
&\e{B}^{(\delta_1,\delta_2)}F(\r{a}[n]; m,\nu)
=\int_0^\infty j_\nu(\delta_1\delta_2mn/u)\cr 
&\times \Big(\int_{\B{R}^2}
F(\r{n}[x_1]\r{w}\r{n}[x_2]\r{a}[u])e(-\delta_1m x_1-
\delta_2nx_2/u)dx_1dx_2\Big) {du\over u}.
}
$$
Thus, 
$$
\eqalignno{ 
&\e{B}^{(\delta_1,\delta_2)}F(\r{a}[n]; m,\nu)
=\int_0^\infty j_\nu(\delta_1\delta_2/u)\cr
&\times \Big(\int_{\B{R}^2}
F(\r{n}[x_1/m]\r{w}\r{n}[mux_2]\r{a}[mnu])e(-\delta_1x_1-
\delta_2x_2)dx_1dx_2\Big)du.\qquad&(6.5.21)
}
$$ 
\medskip
One may desire to compute the double sum $(6.5.16)$ and the last
double integral into closed forms. In the applications 
to ${\cal Z}_2(g)$, we are in a fortuitous situation
that the double sum is transformed into a product of two values of $H_V$. As
to the double integral, it is a Fourier transform over the Euclidean plane, and
thus, in principle, can be expressed in terms of a Bessel transform
as can be seen in $(6.5.24)$ below. With
${\cal Z}_2(g)$, the situation turns out in fact to be as
such. Hence the matter seems to depend much on the specific
nature of the seed $F$.  Nevertheless, with any smooth $F$, one
might appeal to Mellin transform of several variables, and 
$(6.5.21)$ could be pushed to a more closed form. 
\medskip
Before finishing this section,
we render the spectral decomposition of
${\cal Z}_2(g)$ in terms of notions from the theory of
$\varGamma$-automorphic representations:  Thus, let us
put
$$
\eqalignno{
\Theta(\nu;g)&=\int_0^\infty
\Big({u\over u+1}\Big)^{{1\over2}}g_c\left(
\log\left(1+{1/u}\right)\right)
\Xi(u;\nu)d^\times\!u,\qquad&(6.5.23)\cr
\Xi(u;\nu)&=\int_{\B{R}^\times} j_0(-v)j_\nu\left({v\over u}\right)
{d^\times\! v\over\sqrt{|v|}}.&(6.5.24)
}
$$
Then $(4.7.1)$ is expressed as
$$
{\cal Z}_2(g)=\left\{{\cal Z}_2^{(r)}+{\cal Z}_2^{(c)}
+{\cal Z}_2^{(e)}\right\}(g),\eqno(6.5.25)
$$ 
where
$$
\eqalignno{
{\cal Z}_2^{(c)}(g)&=\sum_{V}
\alpha_VH_V(\txt{1\over2})^3\Theta(\nu_V;g),&(6.5.26)\cr
{\cal Z}_2^{(e)}(g)&={1\over 2\pi i}\int_{(0)}
{\left|\zeta\left({1\over2}+\nu\right)\right|^6
\over|\zeta(1+2\nu)|^2}\Theta(\nu;g)d\nu,&(6.5.27) 
}
$$ 
with $\alpha_V=|\varrho_V(1)|^2+|\varrho_V(- 
1)|^2$. The $V$ runs over a complete system of
Hecke-invariant cuspidal irreducible 
$\varGamma$-automorphic representations of
$\r{G}$. The ${\cal Z}_2^{(r)}(g)$ is the same as ${\cal Z}_{2,r}^{(c)}(g)$
in $(4.7.1)$.  The equivalence between 
$(4.7.2)$ and $(6.5.24)$ may independently
be verified by using $(6.4.14)$ and $(6.4.15)$. The factor $j_\nu$ in $(6.5.24)$
has come from the same involved in $(6.5.21)$. 
\vskip 1cm
\centerline{\bf 6.6 Mean values of automorphic $\mib L$-functions}
\medskip
\noindent
The aim of the present section is to develop a unified
treatment of mean values of individual automorphic $L$-functions associated
with the spectral decomposition of $L^2(\varGamma\backslash\r{G})$.
We shall establish a complete spectral expansion for
$$
\e{M}(U,g)=\int_{-\infty}^\infty 
|L_U(\txt{1\over2}+it)|^2g(t)dt,\eqno(6.6.1)
$$
where $U$ is any irreducible representation 
listed in $(6.3.22)$, and 
$$
\matrix{
&\hbox{the weight $g$ is even, entire, real valued on
$\Bbb R$, and}\cr
&\hbox{of rapid decay in any fixed holizontal strip,}
}\eqno(6.6.2)
$$
which is 
assumed for the sake of simplicity and could be replaced by the less stringent
assumption given in the introduction of Chapter 4. 
Our argument is {\it unified\/}
in the sense that it is equally  applicable to any $U$,
whereas hitherto known arguments are applicable only either to the zeta-function,
which corresponds to the continuous spectrum,
or to those $L_U$ with $U$ in the discrete series.
It will be seen that the theory of automorphic representations is 
genuinely needed in our solution of the problem. This is 
in contrast to the situation in the previous section where 
the theory has been utilized to gain a geometric understanding of the sum
formulas and the explicit formula for ${\cal Z}_2(g)$, and 
could be dispensed with otherwise. In other words, the theory of automorphic
representations leads us to a genuinely
new assertion in the theory of mean values
of the zeta and $L$-functions as well. We shall discuss mainly
the case with $U$ in the unitary principal series because of an obvious reason.
\medskip
We begin, nevertheless, with a brief discussion on irreducible
representations in the discrete series, in order to illustrate
the main problematics that we have to resolve when we treat
$\e{M}(U,g)$ with general $U$. Thus, let us assume temporarily that
$U$ be in the holomorphic discrete series and $\nu_U=k-{1\over2}$ with
an integer $k\ge6$. Corresponding to $(4.3.1)$, we consider
$$ 
{\cal J}_U(u,v;g)=\int_{-\infty}^\infty \overline{L_U(\bar{u}+it)}
L_U(v+it)g(t)dt,\eqno(6.6.3)
$$ 
which is an entire function over ${\Bbb C}^2$.
We have $\e{M}(U,g)={\cal J}_U({1\over2},{1\over2};g)$.
In the region of absolute convergence it holds that
$$
{\cal J}_U(u,v;g)={L_{U\otimes U}(2(u+v))\over\zeta(2(u+v))}g^*(0) +
{\cal J}_U^{(1)}(u,v;g)+
\overline{{\cal J}_U^{(1)}(\bar{v},\bar{u};g)},\eqno(6.6.4)
$$  
where
$$  
{\cal J}_U^{(1)}(u,v;g)=\sum_{f,\,n=1}^\infty{{\overline{\varrho_U(n)}
{\varrho_U(n+m)}}\over {n^u (n+m)^v}} g^*\Big(\log {{n+m}\over n}
\Big),\eqno(6.6.5)
$$ 
with $g^*$ as in $(4.1.6)$. By the Mellin inversion, 
$$
\eqalignno{
{\cal J}_U^{(1)}(u,v;g)={1\over{2\pi
i}}\int_{(\eta)}\Big\{ &\sum_{m=1}^\infty m^{-s}D_U(u+v-s,m)\Big \}\cr
&\times\tilde g(s,s-u-k+\txt{3\over2})ds,&(6.6.6)
}
$$ 
where $\tilde{g}$ is defined by $(4.1.7)$, and
$$ 
D_U(s,m)=\sum_{n=1}^\infty{{\overline{\varrho_U(n)}
{\varrho_U(n+m)}}\over {(n+m)^s}}\Big({n\over n+m}\Big)^{k-{1\over2}}.
\eqno(6.6.7)
$$   
If $u+v>\eta+{3\over2}>{5\over2}$, then $(6.6.6)$ converges absolutely, for
we have $(6.3.34)$--$(6.3.35)$. On noting that
$$
|\lambda_U^{(k)}(\r{g})|^2={2^{4k}\pi^{2k+1}\over\Gamma(2k)}y^{2k}
\Big|\sum_{n=1}^\infty\varrho_U(n)n^{k-{1\over2}}e(nz)\Big|^2,\quad
z=x+iy,\eqno(6.6.8)
$$
is in $L^2(\varGamma\backslash{\eusm H})$ and that for $\Re s>1$
$$ 
D_U(s,m)=4{(4\pi)^{s-2}\Gamma(2k)\over{\Gamma(s+2k-1)}}
\langle|\lambda_U^{(k)}|^2, P_m (\cdot,\bar{s})\rangle
\eqno(6.6.9) 
$$
with $P_m$ as in $(1.1.4)$, we may compute a decomposition of
$D_U(s,m)$ over the spectrum of the Casimir operator; 
here the inner product is over $\varGamma\backslash{\eusm H}$, i.e.,
$\varGamma\backslash\r{G}/\r{K}$. Thus, by Theorem 6.1
or rather by Theorem 1.1 together with $(2.1.19)$--$(2.1.20)$, we find that
$$
\eqalignno{
&D_U(s,m)={m^{{1\over2}-s}\Gamma(2k)\over\Gamma(s)\Gamma(s+2k-1)}\cr
&\times\Big\{\sum_V{\raise4pt\hbox{${}\!^{\r{u}}$}}{\varrho_V(m)
\over\pi^{{1\over2}-{\nu_V}}\Gamma({1\over2}+\nu_V)}
\Gamma(s-\txt{1\over2}+\nu_V)
\Gamma(s-\txt{1\over2}-\nu_V)\langle{|\lambda_U^{(k)}|^2,\lambda_V^{(0)}
}\rangle\cr
&+{1\over4\pi i}\int_{(0)}{\sigma_{-2\nu}(m)L^*_{U\otimes
U}(\txt{1\over2}+\nu)\over
|\Gamma({1\over2}+\nu)\zeta(1+2\nu)|^2}\Gamma(s-\txt{1\over2}+\nu)
\Gamma(s-\txt{1\over2}-\nu)d\nu\Big\},&(6.6.10)
}
$$
where $\sum^\r{u}$ is as in $(6.5.1)$, and
$$
L^*_{U\otimes U}(s)=(2\pi)^{2(1-s)}\Gamma(2k)^{-1}\Gamma(s)\Gamma(s+2k-1)
L_{U\otimes U}(s)
\eqno(6.6.11)
$$
is the normalized Rankin $L$-function attached to $U$.
\par
The next step is to insert the decomposition $(6.6.10)$ into $(6.6.6)$,
and try  to exchange the order of the sum and the integral. Here we face a
problem about the uniform growth rate of individual terms on the right
side of $(6.6.10)$, anything similar to which 
we have not experienced in dealing with
${\cal Z}_2(g)$. Thus, in general we may expect at most that the factor
$\tilde{g}(s,s-u-k+{3\over2})$ in $(6.6.6)$ 
decays faster than any negative power of
$|s|$ while $\Re s$ is bounded; 
consequently the polynomial growth of the right side of
$(6.6.10)$ is essential for the success of the argument. Note that the same 
about the contribution of the continuous
spectrum is immediate, via the functional equation 
$L^*_{U\otimes U}(s)=L^*_{U\otimes U}(1-s)$.
Hence we need in turn the polynomial
growth of 
$$
|\langle{|\lambda_U^{(k)}|^2,\lambda_V^{(0)}}\rangle|
\exp(\txt{1\over2}\pi|\nu_V|)\eqno(6.6.12)
$$ 
with respect to the parameter $\nu_V$, in view of
the estimation $(2.3.2)$ and Stirling's formula.
As a matter of
fact, an assertion exists that guarantees such a bound for
$(6.6.12)$. However, the reasoning employed there is highly specific to that
$U$ is in the discrete series, and it does not extend to the general
situation where we have an arbitrary irreducible representation in place of $U$.
Because of this, it is useless for us to 
proceed further along the above argument.  
Nevertheless, it might be worth stating
the following analogue  of $(6.5.23)$--$(6.5.24)$: The contribution of the
irreducible cuspidal representation $V$ to $\e{M}(U,g)$,
$\nu_U=k-{1\over2}$ with an integer $k\ge6$, is equal to
$$
(-1)^k{(2\pi)^{2k-1}\Gamma(2k)H_V(\txt{1\over2})
\langle{|\lambda_U^{(k)}|^2,\overline{\varrho_V(1)}\lambda_V^{(0)}}\rangle
\over\cos(\pi\nu)
\Gamma(2k-{1\over2}+\nu_V)\Gamma(2k-{1\over2}-\nu_V)}
{\Theta_k(\nu_V;g)\over\pi^{{1\over2}-\nu_V}\Gamma({1\over2}+\nu_V)},
\eqno(6.6.13)
$$
where
$$
\eqalignno{
\Theta_k(\nu;g)&=\int_0^\infty
\Big({u\over u+1}\Big)^{1-k}g_c\left(
\log\left(1+{1/u}\right)\right)
\Xi_k(u;\nu)d^\times\!u,\qquad&(6.6.14)\cr
\Xi_k(u;\nu)&=\int_{\B{R}^\times} j_{k-{1\over2}}(-v)
j_\nu\left({v\over u}\right)|v|^{k-1}
{d^\times\! v}.&(6.6.15)
}
$$
\medskip
With this, we now turn to the unitary principal series, so that hereafter
$U$ is an arbitrary irreducible representation with $\nu_U\in i{\Bbb R}$.
One may follow the above argument, with necessary changes, up to $(6.6.7)$.
We face, however, a serious obstacle already at $(6.6.8)$; 
that is, this time we have $\lambda_U^{(0)}$, and 
the factor $K_{\nu_U}(2\pi|n|y)e(nx)$ arises in place of the additive character
$e(nz)$. Hence we  are unable to readily attain an expression analogous to
$(6.6.9)$. On the other hand, despite this
difficulty there exists an argument that extends the  polynomial growth of
$(6.6.12)$ to $|\langle{|\lambda_U^{(0)}|^2,\lambda_V^{(0)}}\rangle|
\exp(\txt{1\over2}\pi|\nu_V|)$; but we 
skip the details, since we are about to exhibit an alternative argument
that resolves as well the above difficulty pertaining to the $K$-Bessel factors. 
\medskip
Our discussion depends much on uniform bounds for
$\e{A}\phi_p(\r{a}[y],\nu)$, $\e{A}=\e{A}^+$, such as $(6.4.17)$ and
$(6.4.18)$. In order to make our argument applicable to any
irreducible cuspidal representation, we derive from $(6.4.18)$ a bound
that is somewhat weaker than $(6.4.17)$ but still sufficient for our
purpose; in fact, the proof of $(6.4.18)$ works for all cases. Let us
assume that $\nu\in i\B{R}$. We divide the integral $(6.4.10)$ at $y
=|p|+|\nu|+1$. To the part with smaller argument we apply the fact that
$\e{A}^{\sgn(u)}\phi_p(\r{a}[|u|];\nu)$ is a unit vector in
$L^2(\B{R}\!^\times\!, d^\times\!/\pi\!)$. Hence this part is
$\ll (|p|+|\nu|+1)^{\Re s-{1\over2}}$. On the other  hand, by 
$(6.4.18)$ the remaining part is $\ll(|p|+|\nu|+1)^{\Re s}$.
We then invoke the identity
$$
\Gamma_p(s,\nu)=4\pi\cdot{\pi\Gamma_p(s+2,\nu)
-p\Gamma_p(s+1,\nu)\over
s^2-\nu^2},\eqno(6.6.16)
$$ 
which can be proved via integration by parts, on noting that
$$
\eqalignno{
&\e{D}_\nu\e{A}\phi_p(\r{a}[y],\nu)=-4\pi p\e{A}\phi_p(\r{a}[y],\nu),\cr
&\e{D}_\nu=(d/dy)^2-(2\pi)^2-\left(\nu^2-\txt{1\over4}\right)\!y^{-2},
&(6.6.17)
}
$$ 
as $\e{A}\phi_p(\r{a}[y],\nu)$ is a constant
multiple of the Whittaker function $W_{p,\nu}(4\pi y)$; see $(3.2.32)$. 
Then, by the Mellin inversion, we have
$$
\e{A}\phi_p(y,\nu)=-2i\int_{(\varepsilon)}
\left(\pi\Gamma_p(s+2,\nu)
-p\Gamma_p(s+1,\nu)\right){y^{{1\over2}-s}\over s^2-\nu^2}ds,\eqno(6.6.18)
$$
with any small constant $\varepsilon>0$.
Inserting the above bound for $\Gamma_p(s,\nu)$, we conclude that
$$
\e{A}\phi_p(\r{a}[y],\nu)
\ll y^{{1\over2}-\varepsilon}(|p|+|\nu|+1)^{2+\varepsilon},\eqno(6.6.19)
$$
where the implied constant depends only on $\varepsilon$.
\par
The discussion on the discrete series is analogous. Actually, any combination of
$p,\nu$ such that either $-{1\over2}<\nu<{1\over2}$ or $\nu=\ell-{1\over2}$ with
$1\le\ell\in\B{Z}$, $\ell\le |p|$,
could also be dealt with, as an explicit evaluation of the norm of 
$\e{A}^{\sgn(u)}\phi_p(\r{a}[|u|];\nu)$ in
$L^2(\B{R}\!^\times\!, d^\times\!/\pi\!)$ can be performed by using
$(6.4.23)$. In passing, we remark that the operator 
$\e{D}_\nu$ is connected with $\partial_\theta$ via the Kirillov map. 
\medskip
We return to $(6.6.3)$ but with the present choice of $U$;
we shall mostly omit the symbol $U$ to avoid otherwise heavy notation, so that
hereafter we have, for instance, $\varrho(n)=\varrho_U(n)$. We then replace
$(6.6.4)$ by
$$
{\cal J}(u,v;g)={L_{U\otimes U}(u+v)\over\zeta(2(u+v))}g^*(0)+J(u,v;g)
+\overline{J(\bar{v},\bar{u};g)},\eqno(6.6.20)
$$
where 
$$
\eqalignno{
&J(u,v;g)\cr
&=\sum_{f=1}^\infty\sum_{n=1}^\infty
{\overline{\varrho(n)}\varrho(n+m)\over (2n+m)^{u+v}}
\left({\sqrt{n(n+m)}\over 2n+m}\right)^{2\alpha}g_*(m/(2n+m);u,v), 
\qquad&(6.6.21)
}
$$
with
$$
g_*(x;u,v)
=2^{u+v+2\alpha}{g^*(\log( (1+x)/(1-x)))\over (1-x)^{u+\alpha}
(1+x)^{v+\alpha}},\quad 0\le x\le1.\eqno(6.6.22)
$$
Here $\alpha$ is a sufficiently large positive integer,
which is implicit throughout the sequel. Let $\tilde{g}$ be the Mellin transform
of $g_*$; note that the definition of $\tilde{g}$ has been
changed from $(4.1.7)$.  It is immediate to see
that $\tilde{g}(s;u,v)$ is of rapid decay with respect to
$s$, provided $\Re s$ and $u,v$ are bounded; moreover,
$\tilde{g}(s;u,v)/\Gamma(s)$ is entire over $\B{C}^3$.  Thus, by Mellin's
inversion, 
$$
J(u,v;g)
={1\over2\pi i}\int_{(\eta)}\!\Big\{\sum_{m=1}^\infty m^{-s}
D(u+v-s,m)\Big\}\tilde{g}(s;u,v)ds,\eqno(6.6.23)
$$
where
$$
D(s,m)=\sum_{n=1}^\infty
{\overline{\varrho(n)}\varrho(n+m)\over (2n+m)^s}
\!\left({\sqrt{n(n+m)}\over 2n+m}\right)^{2\alpha}.\eqno(6.6.24)
$$
It is assumed temporarily that $\Re(u+v)>\max\{2,1+\eta\}$ is 
sufficiently large.
\medskip
We now try to {\it imitate\/} $(6.6.8)$  with a vector 
in $U$ that is generated by $\lambda_U^{(0)}(\r{g})$. What is essential for our
purpose is the fact that the Fourier coefficients $\varrho(n)$ are stable in this
generating process, and the subspace $U$ thus obtained is unitarily equivalent
to the space $L^2({\Bbb R}^\times,\pi^{-1}d^\times)$ as is
stated in Lemma 6.2. With this in mind, 
we apply the inverse Kirillov map $\e{K}^{-1}$ to the function 
$$
w(y,\tau)=\cases{\hfil 0 & if $y\le0$,\cr
y^{\alpha+{1\over2}}\exp(-\tau y) & if $y>0$,}\eqno(6.6.25)
$$
with $\Re\tau>0$,  which is in $L^2({\Bbb R}^\times,\pi^{-1}d^\times)$; all
implicit constants in the sequel may depend on $U$, $\alpha$, and $\Re\tau$ at
most.  Namely, according to the mechanism explained around $(6.4.34)$, we have that
$$
\Phi(\r{g},\tau)=\sum_{\scr{n=-\infty}\atop\scr{n\ne0}}^\infty
{\varrho(n)\over\sqrt{|n|}}\e{A}^{\sgn(n)}\e{K}^{-1}w(\r{a}[|n|]\r{g},\tau)
\eqno(6.6.26)
$$
is a vector in $U$ such that
$$
\Phi(\r{n}[x]\r{a}[y],\tau)=
\sum_{n=1}^\infty{\varrho(n)\over\sqrt{n}}w(ny,\tau)
\exp(2\pi inx).\eqno(6.6.27)
$$
More precisely, we have, by Lemma 6.2,
$$
\Phi(\r{g},\tau)=\sum_{p=-\infty}^\infty a_p(\tau)\lambda_U^{(p)}(\r{g}),
\eqno(6.6.28)
$$
with
$$
a_p(\tau)={1\over\pi}\int_0^\infty w(y,\tau)
\overline{\e{A}\phi_p(\r{a}[y];\nu_U)}\,{dy\over y}.\eqno(6.6.29)
$$
The function $\Phi(\r{g},\tau)$ is regular for $\Re \tau>0$, 
provided $\alpha>2$, since we have
$$
\eqalignno{
\lambda_U^{(p)}(\r{g})&\ll (|p|+1)^2,&(6.6.30)\cr
a_p(\tau)&\ll (|\tau|+1)^{2\alpha}(|p|+1)^{-\alpha}.&(6.6.31)
}
$$
The former can be shown by $(6.4.18)$. 
To prove the latter, we use the operator $\e{D}_{\nu}$ defined in $(6.6.17)$:  
We may assume that $p\ne0$; then,
$$
a_p(\tau)=-{1\over4\pi p}\int_0^\infty w(y,\tau)
\overline{\e{D}_{\nu_U}\e{A}\phi_p(\r{a}[y];\nu_U)}\,dy.\eqno(6.6.32)
$$
We integrate by parts, repeat the procedure $\alpha$ times, and use
the fact that $\Vert\e{A}\phi_p(\r{a}[y];\nu_U)\Vert\le1$ in 
$L^2({\Bbb R}^\times,\pi^{-1}d^\times)$.
\medskip
Next, we put $\Psi(\r{g},\tau)=\Phi(\r{g},\tau)
\overline{\Phi(\r{g},\overline{\tau})}$.
The Parseval formula in $L^2(\varGamma\backslash\r{G})$ gives that
$$
\eqalignno{
\Psi(\r{g},\tau)={3\over\pi}\langle\Psi,1\rangle
&+\sum_V\sum_{p=-\infty}^\infty
\langle\Psi,
\lambda_V^{(p)}\rangle\lambda_V^{(p)}(\r{g})\cr
&+\sum_{p=-\infty}^\infty\int_{(0)}\langle\Psi,
E_p(\cdot,\nu)\rangle E_p(\r{g},\nu){d\nu\over4\pi i}\,,&(6.6.33)
}
$$
with $V$ running over all irreducible cuspidal representations;
the sums over the discrete series need to be
modified appropriately. The convergence is absolute and fast, provided
$\alpha$ is sufficiently large. In fact, we have
$$
\langle\Psi,\lambda_V^{(p)}\rangle\ll (|\tau|+1)^{4\alpha}
(|\nu_V|+|p|)^{-{1\over2}\alpha},\eqno(6.6.34)
$$ 
with which and $(6.6.30)$ the assertion follows. The proof of
this bound and the discussion on the continuous spectrum are to be given
later.
\par
Picking up the $m$-th Fourier coefficient on both sides of $(6.6.33)$
with $\r{g}=\r{n}[x]\r{a}[y]$, while invoking $(6.3.27)$--$(6.3.29)$,
we get, for any $m>0$,
$$
\eqalignno{
y^{2\alpha+1}&\sum_{n=1}^\infty\overline{\varrho(n)}
\varrho(n+m)(n(n+m))^\alpha\exp(-(2n+m)\tau y)\cr
=&\sum_V{\varrho_V(m)\over\sqrt{m}}X_V(my;\tau)
+\int_{(0)}{m^{-\nu}\sigma_{2\nu}(m)\over
\sqrt{m}\zeta(1+2\nu)}Y_\nu(my;\tau){d\nu\over4\pi i}\,,\qquad&(6.6.35)
}
$$
with
$$
\eqalignno{
X_V(y;\tau)&=\sum_{p=-\infty}^\infty\langle\Psi,
\lambda_V^{(p)}\rangle\e{A}\phi_p(\r{a}[y];\nu_V),&(6.6.36)\cr
Y_\nu(y;\tau)&=\sum_{p=-\infty}^\infty\langle\Psi,
E_p(\cdot,\nu)\rangle\e{A}\phi_p(\r{a}[y];\nu).&(6.6.37)
}
$$
A combination of $(6.4.18)$, $(6.6.19)$, and $(6.6.34)$ yields the
uniform bound
$$
X_V(y;\tau)\ll (|\tau|+1)^{4\alpha}
(|\nu_V|+1)^{-{1\over4}\alpha}
y^{{1\over2}-\varepsilon}(y+1)^{-{1\over5}\alpha}.\eqno(6.6.38)
$$
It should be noted that this bound holds for any $V$, since $(6.4.18)$
and $(6.6.19)$ holds for all relevant $\nu$. The
function $Y_\nu$ will be treated later.
\medskip
We are about to verify $(6.6.34)$. We have, by
$(6.3.25)$,
$$
\eqalignno{
&|(\nu_V^2-\txt{1\over4}+i(2p)^2)^q|
|\langle\Psi,\lambda_V^{(p)}\rangle|
=|\langle\Psi,
(\Omega-i\partial_\theta^2)^q\lambda_V^{(p)}\rangle|\cr
=&|\langle(\Omega+i\partial_\theta^2)^q
\Psi,\lambda_V^{(p)}\rangle|
\le\Vert(\Omega+i\partial_\theta^2)^q
\Psi\Vert,&(6.6.39)
}
$$
for any integer $q\ge0$. 
By $(6.1.11)$ and $(6.1.13)$,
$\Omega\lambda_U^{(k)}\overline{\lambda_U^{(l)}}$ is a linear
combination of $\lambda_U^{(k+j)}\overline{\lambda_U^{(l+j)}}$, $j=-1,0,1$,
the coefficients of which are polynomials of the second degree
on $k,l$. Thus, by $(6.6.28)$, 
$$
(\Omega+i\partial_\theta^2)^q\Psi(\r{g})=
\sum_{k=-\infty}^\infty\sum_{l=-\infty}^\infty
b_{k,l}^{(q)}(\tau)\lambda_U^{(k)}(\r{g})
\overline{\lambda_U^{(l)}(\r{g})},\eqno(6.6.40)
$$
where
$$
b_{k,l}^{(q)}(\tau)=\sum_{j=-q}^q d_j(k,l)a_{k+j}(\tau)
\overline{a_{l+j}(\bar\tau)},\eqno(6.6.41)
$$
with a polynomial $d_j(k,l)$ of degree $2q$ on $k,l$; and by $(6.6.31)$,
we have
$$
b_{k,l}^{(q)}(\tau)\ll (|k|+|l|+1)^{2q}{(|\tau|+1)^{4\alpha}
\over((|k|+1)(|l|+1))^{\alpha}},\eqno(6.6.42)
$$
uniformly for $k,l$, and $\tau$ with $\Re\tau>0$. We 
put $q=\!\left[{1\over3}\alpha\right]$, and get the uniform bound
$(\Omega+i\partial_\theta^2)^q\Psi(\r{g},\tau)
\ll(|\tau|+1)^{4\alpha}$, which and $(6.6.39)$ give $(6.6.34)$. 
\medskip
We turn to $Y_\nu$ defined by $(6.6.37)$. We first invoke the
functional equation for $E_p$ that follows from $(3.2.28)$ via
$(6.2.21)$, and  have, for $\Re{\nu}=0$,
$$
\eqalignno{
\langle(\Omega+i\partial_\theta^2)^q\Psi(\cdot,\tau),
&E_p(\cdot,\nu)\rangle=
\pi^{-2\nu}{\zeta(1+2\nu)\over\zeta(1-2\nu)}
{\Gamma\!\left(\txt{1\over2}+\nu+p\right)\over
\Gamma\!\left(\txt{1\over2}-\nu+p\right)}\cr
&\times\int_{\varGamma\backslash\r{G}}
(\Omega+i\partial_\theta^2)^q\Psi(\r{g},\tau)E_{-p}(\r{g},\nu)d\r{g}.
&(6.6.43)
}
$$
Assuming temporarily that $\Re\nu>{1\over2}$, we unfold the last 
integral, and see via $(6.6.40)$ that it is equal to
$$
{L_{U\otimes U}(\nu+\txt{1\over2})\over\zeta(2\nu+1)}
W(\nu,\tau;p,q),\eqno(6.6.44)
$$
where
$$
\eqalignno{
&W(\nu,\tau;p,q)=
\sum_{l=-\infty}^\infty b_{l+p,l}^{(q)}(\tau)\cr
&\times\sum_{\delta=\pm}\int_0^\infty\e{A}^\delta\phi_{l+p}(\r{a}[y],
\nu_U)\overline{\e{A}^\delta\phi_l(\r{a}[y],\nu_U)}
y^{\nu-{3\over2}}dy.&(6.6.45)
}
$$
Again by $(6.4.18)$, $(6.6.19)$, and $(6.6.43)$, we see that 
$W(\nu,\tau;p,q)$ is regular and $\ll(|\tau|+1)^{4\alpha}$ 
for $\Re\tau>0$ and $\Re\nu>-{1\over2}$. Hence, in the same domain,
$$
Y_\nu(y,\tau)=
{L_{U\otimes U}(\nu+\txt{1\over2})\over\zeta(1-2\nu)}
Y^*_\nu(y,\tau),\eqno(6.6.46)
$$
with
$$
Y^*_\nu(y,\tau)=
\pi^{-2\nu}\sum_{p=-\infty}^\infty\!
{W(\nu,\tau;p,q)\over
\left(\nu^2-{1\over4}-i(2p)^2\right)^{q}}
{\Gamma\!\left(\txt{1\over2}+\nu+p\right)\over
\Gamma\!\left(\txt{1\over2}-\nu+p\right)}\e{A}\phi_p(\r{a}[y];\nu).
\eqno(6.6.47)
$$
One may conclude, via $(6.4.17)$--$(6.4.18)$, that
$$
Y_\nu^*(y,\tau)\ll(|\tau|+1)^{4\alpha}
(|\nu|+1)^{-{1\over4}\alpha}
y^{{1\over2}-|\Re\nu|-\varepsilon}(y+1)^{-{1\over5}\alpha},\eqno(6.6.48)
$$
for $\Re\tau>0$ and $\Re\nu>-{1\over2}$.
\medskip
Now we set $\tau=s$. We multiply both sides of $(6.6.35)$ by $y^{s-2}$ and
integrate. We have
$$
\eqalignno{
&D(s,m)=\sum_V m^{{1\over2}-s}\varrho_V(m) \Xi_V(s)\cr
&+\int_{(0)}{m^{{1\over2}-\nu-s}\sigma_{2\nu}(m)\over
\zeta(1+2\nu)\zeta(1-2\nu)}
L_{U\otimes U}\!\left(\nu+\txt{1\over2}\right)
\Upsilon_\nu(s){d\nu\over4\pi i},&(6.6.49)
}
$$
where
$$
\eqalignno{
\Xi_V(s)&={s^{s+2\alpha}\over\Gamma(s+2\alpha)}
\int_0^\infty y^{s-2}X_V(y,s)dy,\cr
\Upsilon_\nu(s)&={s^{s+2\alpha}\over\Gamma(s+2\alpha)}
\int_0^\infty y^{s-2}Y^*_\nu(y,s)dy.&(6.6.50)
}
$$
The bound $(6.6.38)$ implies that $\Xi_V(s)$ is regular and 
$\ll|s|^{4\alpha+{1\over2}}(|\nu_V|+1)^{-{1\over4}\alpha}$
for $\Re s>{1\over2}$. Also, $(6.6.48)$ implies that 
$\Upsilon_\nu(s)$ is regular
and $\ll|s|^{4\alpha+{1\over2}}(|\nu|+1)^{-{1\over4}\alpha}$
for $\Re s>{1\over2}+|\Re\nu|$ and $\Re\nu>-{1\over2}$. 
\medskip
Therefore we have proved 
\medskip
\noindent
{\bf Lemma 6.4.} {\it The function $D(s,m)$ admits the spectral
decomposition $(6.6.49)$ which converges absolutely and
uniformly for $\Re s>{1\over2}$. In particular, $D(s,m)$
is regular and of polynomial growth for $\Re s>{1\over2}$.
} 
\medskip
With this, we return to the expression
$(6.6.23)$ of the function $J(u,v;g)$;
thus we impose $\Re(u+v)>\max\{2,1+\eta\}$ initially. In view
of the fast decay of $\tilde{g}(s;u,v)$, the last
lemma yields immediately that 
$$
\eqalignno{
J(u,v;g)&=\sum_VL_V\!\left(u+v-\txt{1\over2}\right)\Theta_V(u,v;g)\cr
+{1\over4\pi i}&\int_{(0)}{\zeta(u+v-{1\over2}+\nu)
\zeta(u+v-{1\over2}-\nu)\over\zeta(1+2\nu)\zeta(1-2\nu)}\cr
&\times L_{U\otimes U}\left(\txt{1\over2}+\nu\right)
\Lambda_\nu(u,v;g)d\nu, \qquad &(6.6.51)
}
$$
where
$$
\eqalignno{
\Theta_V(u,v;g)&={1\over2\pi i}\int_{(\eta)}
\Xi_V(u+v-\xi)\tilde{g}(\xi;u,v)d\xi,&(6.6.52)\cr
\Lambda_\nu(u,v;g)&={1\over2\pi i}\int_{(\eta)}\Upsilon_\nu(u+v-\xi)
\tilde{g}(\xi;u,v)d\xi. &(6.6.53)
}
$$
We have
$$
\Theta_V(u,v;g)\ll (|\nu_V|+1)^{-{1\over4}\alpha},\quad
\Lambda_\nu(u,v;g)\ll (|\nu|+1)^{-{1\over4}\alpha},\eqno(6.6.54)
$$
where $u,v$ are bounded; the first holds uniformly in the domain
$\Re(u+v)>{1\over2}$, and the second in $\Re(u+v)>{1\over2}+|\Re\nu|$,
$\Re\nu>-{1\over2}$.
\par
We shall discuss the analytic continuation of the expansion $(6.6.51)$.
Let $c>0$ be a sufficiently small constant. We  may move the  
contour in $(6.6.52)$ to $(c)$, provided
$\Re(u+v)>{2\over3}$, say. Hence, the expansion $(6.6.51)$ holds under 
$\Re (u+v)>{3\over2}$. This condition
is required to get the factors $L_V\!\left(u+v-{1\over2}\right)$ and
$\zeta\!\left(u+v-{1\over2}\pm\nu\right)$. However,
the former is  entire and of a polynomial order
in $\nu_V$ if $u,v$ are bounded. Thus the cuspidal part of
$J(u,v;g)$ is regular in a neighbourhood of 
the point $\left({1\over2},{1\over2}\right)$ at which
it takes the value
$$
\sum_V L_V\!\left(\txt{1\over2}\right)\Theta_V(g),\eqno(6.6.55)
$$
with $\Theta_V(g)=\Theta_V\left({1\over2},{1\over2};g\right)$.
\par
As we are about to deal with the continuous spectrum, 
we should remark that $\Lambda_\nu(u,v;g)$ remains 
regular in the three complex variables  and 
of fast decay in $\nu$, throughout the
procedure below, because of the property of 
$\Upsilon_\nu(s)$ mentioned above. Thus, we temporarily
restrict $(u,v)$ so that $2>\Re(u+v)>{3\over2}$.
Then, in $(6.6.51)$ one may shift the contour to
$({1\over2}+c)$, with $c$ as above, encountering the pole at
$u+v-{3\over2}$ with the residue 
$$
-{L_{U\otimes U}(u+v-1)\over\zeta(2(2-u-v))}
\Lambda_{u+v-\txt{3\over2}}\!\left(u,v;g\right)\eqno(6.6.56)
$$
as well as those of the factor 
$L_{U\otimes U}({1\over2}+\nu)/\zeta(1-2\nu)$; 
we may assume, without any loss of generality, that
$u,v$ are such that all the residues in question are finite.
This yields a meromorphic continuation of the continuous spectrum
part,  so that one may move
$(u,v)$ close to $\left({1\over2},{1\over2}\right)$ as far as
$\Re(u+v)>1$ is satisfied, which 
is needed to have the last $\Lambda$ factor
defined well.  Then, shift the
$\nu$-contour back to the original. All the residual contribution
coming from $L_{U\otimes U}({1\over2}+\nu)/\zeta(1-2\nu)$ cancel out those
arising from the previous shift of the contour. Only the
pole at ${3\over2}-u-v$ contributes newly. The resulting
integral is regular at $\left({1\over2},{1\over2}\right)$; we get the
factor $\Lambda_\nu(g)=\Lambda_\nu\!\left({1\over2},
{1\over2};g\right)$, $\Re\nu=0$.
\medskip
\medskip
Collecting all the above, we obtain
\medskip
\medskip
\noindent
{\bf Theorem 6.3} {\it  
We have the spectral decomposition
$$
\eqalignno{
&\e{M}(U, g)=m(U,g)+2\Re\Bigg\{\!\sum_V L_V
\!\left(\txt{1\over2}\right)\Theta_V(g)\cr
+&\int_{(0)}\!{\zeta\!\left({1\over2}+\nu\right)\!
\zeta\!\left({1\over2}-\nu\right)\over\zeta(1+2\nu)\zeta(1-2\nu)
}L_{U\otimes U}\!\left(\txt{1\over2}+\nu\right)\!\Lambda_\nu(g)
{d\nu\over4\pi i}\!\Bigg\},&(6.6.57)
}
$$
where 
$m(U,g)$ is the value at $\left({1\over2},{1\over2}\right)$ of
the function
$$
\eqalignno{
&{L_{U\otimes U}(2(u+v))\over\zeta(2(u+v))}
g^*(0)\cr
+\Big\{&L_{U\otimes U}(u+v-1)\Lambda_{u+v-\txt{3\over2}}\!
\left(u,v;g\right)\cr
+&L_{U\otimes U}(1-u-v)\Lambda_{\txt{3\over2}-u-v}\!\left(u,v;g\right)\cr
+&L_{U\otimes U}(u+v-1)\Lambda_{u+
v-\txt{3\over2}}\!\left(v,u;g\right)\cr
+&L_{U\otimes
U}(1-u-v)\Lambda_{\txt{3\over2}-u-v}
\!\left(v,u;g\right)\!\Big\}/\zeta(2(2-u-v)). &(6.6.58)
}
$$
}
\bigskip
\noindent
Albeit the $\Lambda$ factors in the last expression is defined so far
only under the condition $2>\Re(u+v)>1$, the expression can
in fact be continued to a neighbourhood of
$\left({1\over2},{1\over2}\right)$, for ${\cal J}(u,v;g)$ and all other parts
in the spectral expansion of $J(u,v;g)$ and 
$\overline{J(\bar{v},\bar{u};g)}$ are regular
there. 

\bye